\documentclass[11pt]{article}
\usepackage{mathrsfs}
\usepackage{latexsym,amsfonts,amssymb,amsmath,amsthm}
\usepackage{graphicx,color}
\usepackage[normalem]{ulem}
\allowdisplaybreaks[4]
\usepackage{array, makecell,caption, threeparttable,boxedminipage}
\usepackage{subeqnarray}
\usepackage{cases}
\usepackage{multirow}
\usepackage{indentfirst}
\usepackage{hyperref}
\usepackage[ruled,vlined]{algorithm2e}
\hypersetup{
	colorlinks = true,
	linkcolor = blue,
	anchorcolor =blue,
	citecolor = blue,
	filecolor = blue,urlcolor = black,
	pdfauthor=author
}
\usepackage{cleveref}

\newcommand{\figref}[1]{Fig.~\ref{#1}}

\newcommand{\tabref}[1]{Tabel~\ref{#1}}

\renewcommand{\theequation}{\thesection.\arabic{equation}}

\begin{document}
	\setlength{\baselineskip}{16pt}

	\setlength{\parindent}{2em}
	\newtheorem{theorem}{Theorem}[section]
	\newtheorem{remark}{Remark}[section]
	\newtheorem{proposition}{Proposition}[section]
	\newtheorem{lemma}{Lemma}[section]
	\newtheorem{definition}{Definition}[section]
	\renewcommand{\theequation}{\thesection.\arabic{equation}}
	\renewcommand{\thefootnote}{}
	\numberwithin{equation}{section}
	\renewcommand{\labelenumi}{(\arabic{enumi})}

	\newcommand\del[1]{}
	
	\title{\Large\bf Randomized Radial Basis Function Neural Network for Solving  Multiscale Elliptic Equations  
	}
	
	\author{Yuhang Wu$^{a,b}$, Ziyuan Liu$^a$, Wenjun Sun$^b$, Xu Qian$^{a,}$\thanks{Corresponding author.}\\
		\\
		{\small  $^a$ \emph{College of Sciences, National University of Defense  Technology}, }\\
		{\small  \emph{Changsha, Hunan 410073, China}}\\
		{\small $^b$ \emph{Institute of Applied Physics and Computational Mathematics,} }\\
		{\small \emph{Beijing 100094, China}}\\
	}
	\date{}
	\maketitle
	\let\thefootnote\relax\footnotetext{\emph{Email addresses}: \href{wuyuhang18@nudt.edu.cn}{\textrm{\bf{wuyuhang18@nudt.edu.cn}}} (Yuhang Wu), \href{liuziyuan17@nudt.edu.cn}{\textrm{\bf{liuziyuan17@nudt.edu.cn}}} (Ziyuan Liu), \href{sun\_wenjun@iapcm.ac.cn}{\textrm{\bf{sun\_wenjun@iapcm.ac.cn}}} (Wenjun Sun),  \href{qianxu@nudt.edu.cn}{\textrm{\bf{qianxu@nudt.edu.cn}}} (Xu Qian)}
	\noindent\hrulefill

		\noindent\textbf{Abstract:} Ordinary deep neural network-based methods frequently encounter difficulties when tackling multiscale and high-frequency partial differential equations. To overcome these obstacles and improve computational accuracy and efficiency, this paper presents the Randomized Radial Basis Function Neural Network (RRNN), an innovative approach explicitly crafted for solving multiscale elliptic equations. The RRNN method commences by decomposing the computational domain into non-overlapping subdomains. Within each subdomain, the solution to the localized subproblem is approximated by a randomized radial basis function neural network with a Gaussian kernel. This network is distinguished by the random assignment of width and center coefficients for its activation functions, thereby rendering the training process focused solely on determining the weight coefficients of the output layer. For each subproblem, similar to the Petrov-Galerkin finite element method, a linear system will be formulated on the foundation of a weak formulation. Subsequently, a selection of collocation points is stochastically sampled at the boundaries of the subdomain, ensuring satisfying $C^0$ and $C^1$ continuity and boundary conditions to couple these localized solutions. The network is ultimately trained using the least squares method to ascertain the output layer weights. To validate the RRNN method's effectiveness, an extensive array of numerical experiments has been executed. The RRNN is firstly compared with a variety of deep neural network methods based on gradient descent optimization. The comparative analysis demonstrates the RRNN's superior performance with respect to computational accuracy and training time. Furthermore, it is contrasted with the Local Extreme Learning Machine method, which also utilizes domain decomposition and the least squares method. The comparative findings suggest that the RRNN method can attain enhanced accuracy at a comparable computational cost, particularly pronounced in scenarios with a smaller scale ratio $\varepsilon$.

	\vskip 2mm
	\noindent\emph{Keywords: } Multiscale elliptic equations; Randomized radial basis function neural network; Domain decomposition; Weak formulation; Least squares method.
	
	\noindent\hrulefill
	

	\linespread{1.2}
	\section{Introduction}

	In the realms of scientific and engineering research, phenomena often exhibit characteristics across multiple scales, presenting a significant challenge in the form of multiscale problems. These problems are frequently modeled as partial differential equations (PDEs) with coefficients that exhibit rapid oscillations and are common in many disciplines, such as chemistry\cite{Bolnykhchemistry}, biology\cite{Dadabiological,sunbiological}, materials science\cite{Linfiber}, and fluid dynamics\cite{Changfluid,Pitschfluid}. For instance, in the context of fluid dynamics, the intricate interplay of vortices across various scales results in the formation of complex flow patterns. Conventional numerical methods face considerable difficulties in accurately simulating these phenomena, primarily due to the necessity of employing an exceedingly fine mesh to capture the high-frequency components of the solution.
	
	In this paper, we will mainly focus on the classical multiscale elliptic problem with Dirichlet boundary conditions:
	\begin{eqnarray}\label{1.1}
		\begin{cases}
			\mathcal{L}u^\varepsilon(\boldsymbol{x}) = f(\boldsymbol{x}), \qquad \boldsymbol{x} \in \Omega,\\
			\quad   u^\varepsilon(\boldsymbol{x}) = g(\boldsymbol{x}), \qquad \boldsymbol{x} \in \partial\Omega,\\
		\end{cases}
	\end{eqnarray}
	where $ \mathcal{L} $ is a linear operator and  $ \mathcal{L}u^\varepsilon(\boldsymbol{x}) = -\text{div}(A^\varepsilon(\boldsymbol{x})\triangledown u^\varepsilon(\boldsymbol{x})) $. $ \Omega  $ is a bounded subset of $ \mathbb{R}^n $ and $ \varepsilon \ll 1   $ is a parameter that represents the ratio of the smallest scale to the largest scale in the physical problems. $ A^\varepsilon(\boldsymbol{x}) = (a^\varepsilon_{ij}(\boldsymbol{x}) ) $ is assumed to be a symmetric, positive definite matrix satisfying
	$$  \lambda|\xi|^2 \le a^\varepsilon_{ij}(\boldsymbol{x}) \xi_i \xi_j  \le \Lambda |\xi|^2 , \qquad \forall \xi \in \mathbb{R}^n,\boldsymbol{x} \in \bar{\Omega}$$
	for some positive constants $\lambda  $ and  $\Lambda $ with $  $$ 0<\lambda <\Lambda $.
	
	This equation has attracted a great deal of attention because of its applicability and simplicity. In recent years, an enormous amount of different conventional numerical methods have been proposed, such as homogenization\cite{Dorobantuhomogenization, FishMultigrid,OwhadiHomogenization} , multiscale finite element method (MsFEM)\cite{Chenmixed, Houmultiscale,HouRemoving}, orthogonal decomposition\cite{PeterseimLocalization}, domain decomposition\cite{Aarnesdomain}, adaptive local bases\cite{WeymuthAdaptive}. When simulating multiscale problems numerically, conventional methods face challenges in accuracy and computational cost due to the generation of complex meshes.
	
	For the past few years, with the successful application of machine learning in science and engineering computing, neural network-based methods have been proposed for solving PDEs. A significant advantage of neural network-based methods is avoiding the complex mesh generation process, which can be challenging for high-frequency problems. These methods train the neural networks to find a optimal solution of PDEs by minimizing the loss function. For instance, the deep Galerkin method (DGM)\cite{dgm} and the physics informed neural network (PINN)\cite{pinns} substituted the neural network into the equation to get the direct residuals of the equation, and minimized the residuals to get a neural network to approximate the true solution. While  the deep Ritz method (DRM)\cite{drm} defined the loss function in energy  formulation and used the principle of minimal potential energy to transform the PDEs solution into an optimization problem.
	
	Despite the significant potential and capability demonstrated by neural network-based methods in solving PDEs, particularly in high-dimensional cases, these approaches often encounter challenges when addressing multiscale  and high-frequency PDEs. Xu\cite{XuTraining,XuFrequency}  investigated the training behavior of deep neural networks (DNNs) in the frequency domain, leading to the proposal of the Frequency Principle (F-Principle). This principle reveals that many DNNs tend to learn the low-frequency components of data with good generalization error, rather than the high-frequency components, which is a primary reason for the low accuracy observed in solving multiscale elliptic problems using conventional DNNs.
	
	To harness the considerable potential of DNNs in scientific and engineering computations involving PDEs, there is substantial value in developing DNN-based methods for multiscale problems. The multiscale deep neural network (MscaleDNN) method\cite{LiDNN,MscaleDNN} has been introduced to mitigate these issues to some extent. It employs radial scaling to transform the original high-frequency data into a low-frequency space and utilizes compact support activation functions to facilitate the separation of frequency content. Building on this, the same team integrated the MscaleDNN method with traditional numerical analysis principles to develop the subspace decomposition-based DNN (SD$^2 $NN)\cite{LiSubspace} for a class of multiscale problems. The SD$^2 $NN incorporates high-frequency MscaleDNN sub-modules and low-frequency standard DNN sub-modules, enabling it to capture both the oscillatory and smooth components of multiscale solutions. In addition to these deep network-based methods, the physics-informed radial basis function network has been developed to maintain local approximation properties throughout the training process. Bai\cite{Bai2023} introduced the physics-informed radial basis network (PIRBN) method, which inherently possesses the local approximation property and has shown superior performance compared to PINN in challenging PDEs with high-frequency features. Addressing multiscale problems, Wang and Chen proposed the sparse radial basis function neural network (SRBFNN)\cite{srbfnn} method, applying a simple three-layer radial basis function neural network to approximate solutions to multiscale problems and significantly accelerating the network's training process.

	Training these neural network based on gradient descent method involves a trade-off between accuracy and efficiency. Increasing the number of training epochs can enhance approximation accuracy but also prolongs training time. Fortunately, extreme learning machines (ELM)\cite{elm1,elm2,elm0} can effectively address some of these issues. ELM employs random weight/bias coefficients and utilizes the least squares method for training, which significantly reduces the number of parameters to be trained and avoids local optima. In recent years, numerous ELM-based methods have been proposed for solving PDEs, achieving high accuracy and substantially reducing training time\cite{rfm,rfmtime,locelm,pielm, elmnet,Shang2023,sunJin}.

	The complexity in solving multiscale problems stems from the distinct local features of the solution, the abrupt gradient changes, and the difficulty in capturing high-frequency information. To more effectively solve the linear multiscale elliptic equations \eqref{1.1} using a neural network, we propose the \textbf{R}andomized \textbf{R}adial basis function \textbf{N}eural \textbf{N}etwork (RRNN) method, drawing on the theory of ELM and the excellent local approximation properties of the radial basis function neural network method\cite{MaiDuyMesh,Moseley2023}. The RRNN method firstly employs domain decomposition, and each subdomain's solution approximated by a local radial basis function neural network that uses Gaussian radial basis functions as activation functions. This approach facilitates better capture of the solution's local properties.  Inspired by ELM, the width and center coefficients of the Gaussian activation function are randomly selected from a uniform distribution. The RRNN method employs a variational formulation in each subdomain. Similar to Petrov-Galerkin finite elements method, it utilizes the random radial basis function neural network space approximation as the trail function space, and can flexibly choose the appropriate test function space. In addition, by applying a set of collocation points within each subdomain boundary and enforcing continuity conditions or boundary conditions on these points,  a final linear system of equations for the output layer's weight parameters is jointly established.  This system is solved using the least squares method, enabling the radial basis function neural network to provide the solution to problems \eqref{1.1}. We present numerical examples to demonstrate the efficacy of the proposed RRNN method. It outperforms existing deep neural network-based methods that rely on gradient descent method, including PINN, DRM, MscaleDNN, and SRBFNN,  in terms of both accuracy and training time reduction. Furthermore, when compared to ELM-based methods like LocELM, the RRNN method achieves lower errors with a similar computational cost, particularly at smaller scale ratio $ \varepsilon$. To our knowledge, this is the first application of ELM-based methods with radial basis function to multiscale elliptic problems. The key innovations of this work are as follows:

	\begin{itemize}
		\setlength{\itemsep}{0pt}
		\setlength{\parsep}{0pt}
		\setlength{\parskip}{0pt}
		\item[(i)] By domain decomposition and normalization, the original high-frequency data can be transformed into a low-frequency space. Then we approximate the sub-solution by using a local radial basis function neural network, which naturally possesses the local approximation property and makes it possible to better capture the local properties of the solution.
		\item[(ii)] The loss function will be constructed by the weak  formulation, and the derivation of $ A^\varepsilon(\boldsymbol{x}) $ and $ \triangledown u^\varepsilon(\boldsymbol{x}) $ will be transferred to a test function $v(\boldsymbol{x})  $, which can avoid the effect of drastically varying gradients on the solution.
		\item[(iii)] The coefficients of the Gaussian activation function are pre-set randomly with a uniform distribution and fixed in the training process, which reduces the total number of training parameters and extremely reduces the training time. In addition, the problem is transformed into a system of linear equations, which can be solved with high accuracy by the linear least squares method. 	
	\end{itemize}	
	
	The remainder of this paper is organized as follows: Section 2 provides a comprehensive review of the radial basis function neural network (RBFNN) framework and reiterates the theoretical underpinnings of the extreme learning machine with RBF kernels.	In Section 3, we elaborate on the framework of the proposed RRNN methodology and detail its implementation.
	Section 4 presents a series of numerical experiments, both one-dimensional and two-dimensional, designed to substantiate the effectiveness and superiority of the RRNN method over several neural network-based methods. Finally, Section 5 offers concluding remarks and discusses potential directions for future research.

	\section{Randomized radial basis function neural network}

	\subsection{Radial basis function neural network}
\begin{figure}[thbp!]
	\centering
	\includegraphics[width=1\linewidth]{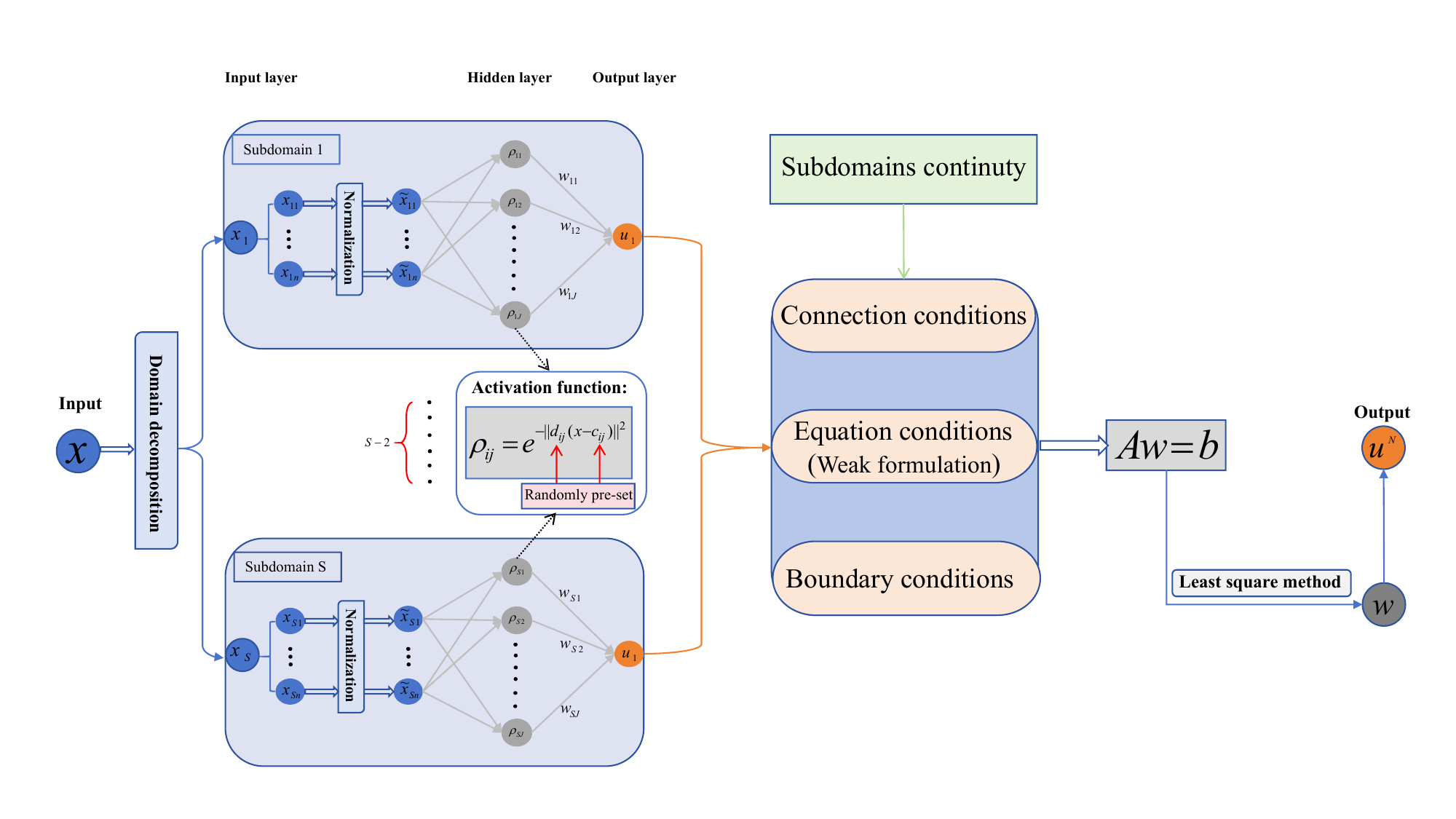} 
	\caption{\footnotesize The structure of randomized radial basis function neural network for solving problems \eqref{1.1}.}
	\label{fig1}
\end{figure}

		A neural networks aims to learn a nonlinear mapping from the input space to the output space. One of the most popular high-performance neural networks is the radial basis function neural network, which is a special type of three-layer artificial neural network. As shown in \figref{fig1}, the layers are input layer, hidden layer and output layer. 
	
 The single hidden layer uses a radial basis function (RBF) as the activation function. RBF is a real-valued function whose value depends only on the distance from the the center point, denoted as $\rho(\boldsymbol{x})= \rho(||\boldsymbol{x}-\boldsymbol{c}||) $ , where $ ||\cdot|| $ usually denotes the standard Euclidean distance.  In this paper, we use the most commonly used radial basis functions: Gaussian function $ \rho(\boldsymbol{x},\boldsymbol{c},\boldsymbol{d})=e^{- \|\boldsymbol{d}(\boldsymbol{x}-\boldsymbol{c})\|^2} .$ where $ \boldsymbol{x}=(x_1,...,x_n)^T $ is the input variable, $\boldsymbol{d}=\text{diag}(d_1,...,d_n)  $ is width coefficient that controls the radial range of action of the function and $ \boldsymbol{c}=(c_1,...,c_n)^T $ is the center coefficient. 
 
 The output value of the activation function (radial basis function) depends on the distance of the input variable from the center. This means that as the distance $ ||\boldsymbol{x}-\boldsymbol{c}|| $ between the point $ \boldsymbol{x} $ and the center $ \boldsymbol{c} $ gets closer to 0, then the contribution of this point gets closer to 1, whereas the contribution gets closer to 0 if the distance $ ||\boldsymbol{x}-\boldsymbol{c}|| $ increases.  In other words, the radial basis function represents the local receptor. For a given input value $ \boldsymbol{x} $, only a small portion near the center is activated. Thus, radial basis function neural network has a good local mapping property.

The linear operation of the output layer is given as follows:
\begin{equation}\label{2.1}
		u^{N}(\boldsymbol{x})=\sum_{i=1}^{J} w_i \rho(\boldsymbol{x},\boldsymbol{c}_i,\boldsymbol{d}_i)=\sum_{i=1}^{J} w_i e^{- \|\boldsymbol{d}_i(\boldsymbol{x}-\boldsymbol{c}_i)\|^2} ,
\end{equation}
where $ u^{N}(\boldsymbol{x}) $ is the output for the input vector $ \boldsymbol{x} $, $ J $ is the number of neurons of hidden layer, and $ \boldsymbol{c}_i,\boldsymbol{d}_i$ are center and width coefficients for the $i$th Gaussian function.

\begin{remark}
	 \rm  In this paper, for convenience, we consider $ d_{i1}=d_{i2}=\cdots=d_{in} $ for the width coefficient $ \boldsymbol{d}_i $ of the $i$th Gaussian function. Then \eqref{2.1} can be rewrote as:
	 \begin{equation}\label{2.2}
	 	u^{N}(\boldsymbol{x})=\sum_{i=1}^{J} w_i e^{-\sigma_i \|\boldsymbol{x}-\boldsymbol{c}_i\|^2} ,
	 \end{equation}
 where $ \sigma_i  $ is called the shape coefficient of $i$th  Gaussian function and $ \sigma_i=d_{in}^2 .$
\end{remark}

Neural networks can be considered as approximation schemes due to universal approximation theorem, and RBFNN also has good approximation properties for continuous functions. Poggio and Girosi\cite{PoggioNetworks} have proved that the multi-layer networks of the back-propagation type do not have the best approximation property. Meanwhile, they proved the existence and uniqueness of the best approximation for RBFNN. Compared with Back-Propagation neural network(BPNN), an important advantage of RBFNN is the existence of a fast linear learning algorithm in a network capable of representing complex nonlinear mappings, greatly speeding up the learning of the neural network while avoiding the local minima problem\cite{Moodylearning}.

\subsection{Randomized radial basis function neural network}

\indent Traditional feed-forward neural networks often necessitate solving a comprehensive nonlinear optimization problem, a process that tends to be slower than desired. This has historically been a limitation in their practical application. The RBFNN addresses this issue to some extent by initially employing an appropriate unsupervised algorithm to select the optimal width, center and weight coefficients of RBFNN, such as  K-means clustering algorithm\cite{Lim2008}, recursive orthogonal least squares\cite{Gomm2000} and the quick propagation method\cite{Montazer2007}. Thereby circumventing the need for a complete nonlinear optimization of the network, which significantly speeds up the training process of the network. 

Despite these improvements, the computational expense associated with RBFNNs remains relatively high. In the past few decades, the introduction of the extreme learning machine (ELM)\cite{elm0} has significantly mitigated this challenge. The ELM is distinguished by its good generalization performance and rapid learning capabilities. It operates by randomly assigning the coefficients of the hidden layer neurons and then analytically determining the output weights for single hidden layer feed-forward neural networks (SLFNs). The universal and approximation capabilities of ELMs with various activation functions have been extensively studied in the works of Huang\cite{elm1,elm2}.
Building upon these advancements, this paper introduces the RRNN method, which utilizes RBF kernels with randomly assigned parameters. The RRNN method is explored as an innovative approach to further enhance the efficiency and accuracy of neural network-based solutions for multiscale PDEs.

In RRNN, the output function with $ J $ hidden neurons can be represented by 
\begin{equation}\label{2.3}
	\boldsymbol{u}^{N}(\boldsymbol{x}) = \sum_{i=1}^{J}\boldsymbol{w}_i \rho(\boldsymbol{c}_i,\sigma_i,\boldsymbol{x}),\qquad \boldsymbol{x}\in \mathbb{R}^n, \boldsymbol{w}_i\in \mathbb{R}^m, \boldsymbol{c}_i\in \mathbb{R}^n, \sigma_i \in \mathbb{R}^+,
\end{equation}	
where $  \rho(\boldsymbol{c}_i,\sigma_i,\boldsymbol{x}) $ denotes the activation function of $i$th hidden neuron.

For $ Q $ arbitrary distinct samples $ (\boldsymbol{x}_i,\boldsymbol{u}_i) \in \mathbb{R}^n \times \mathbb{R}^m $, RRNN with $ J $ hidden neurons are mathematically modeled as 
\begin{equation}\label{2.4}
	\sum_{i=1}^{J}\boldsymbol{w}_i \rho(\boldsymbol{c}_i,\sigma_i,\boldsymbol{x_j}) = \boldsymbol{u}_j^{N}, \qquad j=1,...,Q.
\end{equation}
By the approximation theorem, there exists $\boldsymbol{w}_i, \boldsymbol{c}_i$ and $\sigma_i $ such that RRNN can approximate these $ Q $ samples with zero error means that
\begin{equation}\label{2.5}
	\sum_{i=1}^{J}\boldsymbol{w}_i \rho(\boldsymbol{c}_i,\sigma_i,\boldsymbol{x}_j) = \boldsymbol{u}_j, \qquad j=1,...,Q,
\end{equation}
which can be written compactly as:
\begin{equation}\label{2.6}
	\boldsymbol{A}\boldsymbol{w}=\boldsymbol{b},
\end{equation}
where
\begin{equation}\label{2.7}
	\boldsymbol{A}=\begin{bmatrix}
		\rho(\boldsymbol{c}_1,\sigma_1,\boldsymbol{x}_1)& \cdots & \rho(\boldsymbol{c}_J,\sigma_J,\boldsymbol{x}_1)\\
		\vdots & \cdots  &\vdots  \\
		\rho(\boldsymbol{c}_1,\sigma_1,\boldsymbol{x}_Q)& \cdots & \rho(\boldsymbol{c}_J,\sigma_J,\boldsymbol{x}_Q)
	\end{bmatrix}_{Q\times J,}
\boldsymbol{w}=\begin{bmatrix}
	\boldsymbol{w}^T_1 \\
	\vdots \\
	\boldsymbol{w}^T_J
\end{bmatrix}_{J\times m,}
\boldsymbol{b}=\begin{bmatrix}
	\boldsymbol{u}^T_1 \\
	\vdots \\
	\boldsymbol{u}^T_Q
\end{bmatrix}_{Q\times m.}
\end{equation}

If the activation function $\rho  $ is infinitely differentiable,  the parameters of the  hidden layer output matrix $ \boldsymbol{A} $ can be randomly generated. Huang\cite{elm2} have proved its functional approximation capability, which can be derived  from the following fundamental theorem.

\begin{theorem}
Given any small positive value $ \epsilon>0 $, the activation function $ \rho : \mathbb{R}^n \rightarrow \mathbb{R}  $ which is infinitely differentiable in any interval, and $ Q $ arbitrary distinct samples $ (\boldsymbol{x}_i,\boldsymbol{u}_i) \in \mathbb{R}^n \times \mathbb{R}^m $, there exists $ J \le Q $ such that for any $ \{ \boldsymbol{c}_i,\sigma_i \}_{i=1}^J $ randomly generated from any intervals of $ \mathbb{R}^n \times \mathbb{R} $, according to any continuous probability distribution, with probability one $ \|A\boldsymbol{w}-\boldsymbol{b} \|\le \epsilon. $
\end{theorem}

Theorem 2.1 shows that for any $ \{ \boldsymbol{c}_i,\sigma_i \}_{i=1}^J $ randomly generated from any intervals of $ \mathbb{R}^n \times \mathbb{R} $ with any continuous probability distribution, the RBFNN output \eqref{2.2} has a good approximation capability. Therefore, we can impose further assumptions for the shape and center coefficients for the hidden layer that they are pre-set to uniform random values generated on an interval of $ \mathbb{R}$ and $\mathbb{R}^n$.  More specifically, we assume that:
\begin{equation}\label{2.8}
	\sigma_{i} \in  \mathbb{U}([0,\beta])  ,\qquad \boldsymbol{c}_{i} \in  \mathbb{U}([-1,1]^n), \qquad i=1,...,J,
\end{equation}
where $ \beta $ is a positive constant and referred to as the upper bound of random shape coefficient $ \sigma $.

Once the two coefficients are randomized  with Equ.\eqref{2.8}, they will be fixed throughout the training and cannot be change, and the output becomes a linear combination of several $ \rho(\boldsymbol{x}) $. Then $ \{\rho_i(\boldsymbol{x})\}_{i=1}^{J} $ can be regarded as a group of basis functions. Hence, wo denote the randomized radial basis function neural network space 
\begin{equation}\label{2.9}
	U_R(\Omega)  = \Big\{  u^N(\boldsymbol{w},\boldsymbol{x})  = \sum_{i=1}^{J}\boldsymbol{w}_i \rho_i(\boldsymbol{x}):\boldsymbol{x} \in \Omega  \Big\}= \text{span}\{  \rho_1(\boldsymbol{x}), \rho_2(\boldsymbol{x}),..., \rho_J(\boldsymbol{x})\}  .
\end{equation}

	\section{Method }
	In this section, we will detail the difficulty of solving the multiscale elliptic equations  \eqref{1.1} and how the RRNN method can solve it well.
	
\begin{figure}[htbp!]
		\centering
		\footnotesize
		\begin{tabular}{@{\extracolsep{\fill}}c@{}c@{\extracolsep{\fill}}}
			\includegraphics[width=0.33\linewidth]{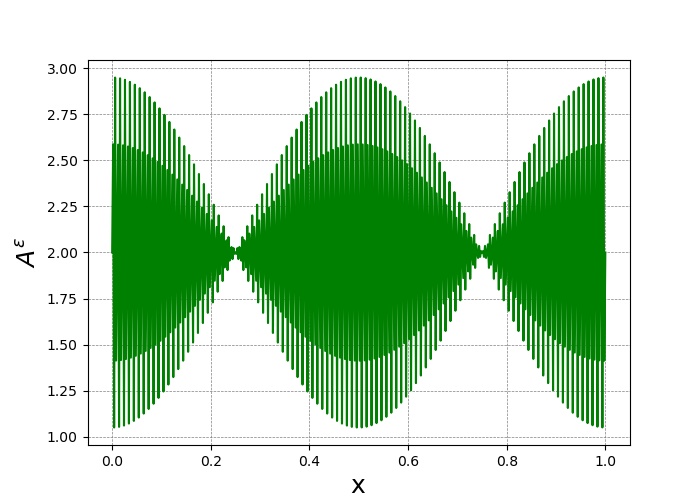} &
			\includegraphics[width=0.33\linewidth]{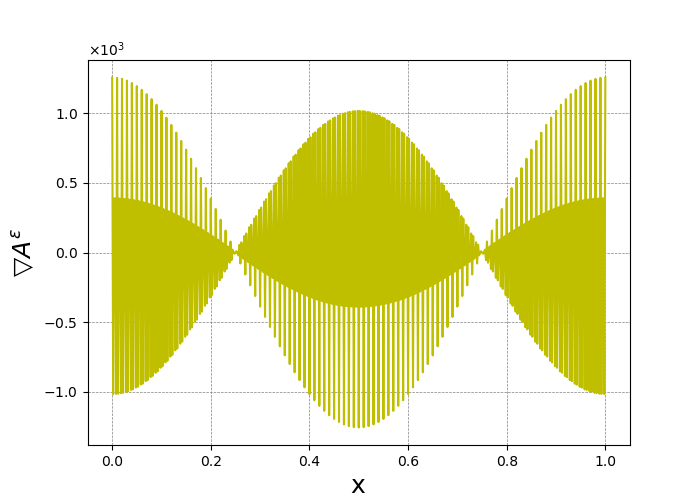}  \\
			(a) coefficient $A^\varepsilon(x)$  & (b) derivative of coefficient  $\triangledown A^\varepsilon(x)$  \\
			\includegraphics[width=0.33\linewidth]{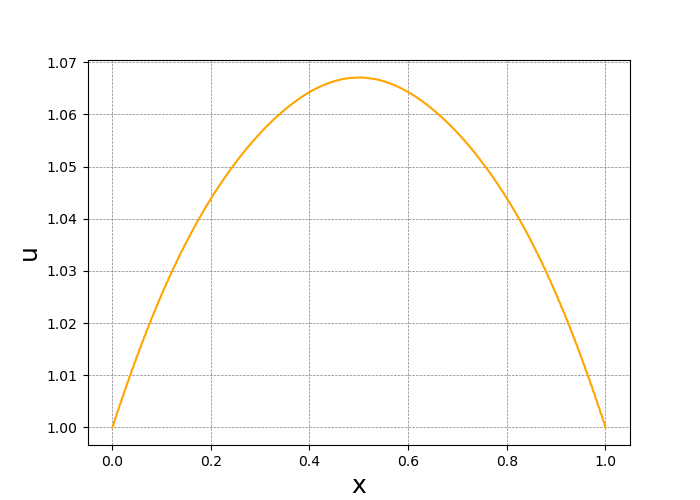} &
			\includegraphics[width=0.33\linewidth]{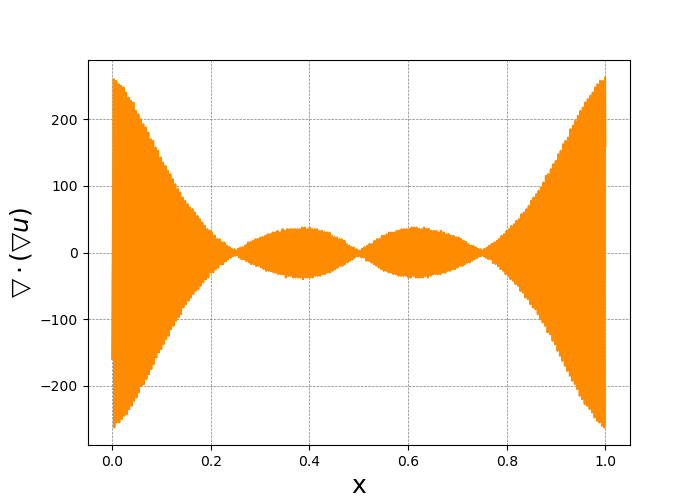}  \\
		(c) solution $u(x)$  & (d) second-order derivative   $\triangledown \cdot (\triangledown u(x))$  \\			
		\end{tabular}
		\caption{\footnotesize The coefficient $A^\varepsilon(x)$ and solution and their derivative of Example 4.2 in \textbf{Section 4} when $\varepsilon=0.005.$ }
		\label{fig21}
	\end{figure}

	\textbf{The difficulty of solving the multiscale elliptic equations.} The multiscale elliptic equations  \eqref{1.1} has a coefficient with a small scale ratio $\varepsilon$, which exhibits violent oscillatory characteristics. Taking Example 4.2 in \textbf{Section 4} as  an example, it can be seen from \figref{fig21}(a) that the coefficient $A^\varepsilon(x)$ oscillates violently in the interval $ x\in[0,1]$ when $\varepsilon$ = 0.005, with strong high-frequency characteristics. In addition, the equation contains the derivative term of $A^\varepsilon(x)$, shown in \figref{fig21}(b), which also oscillates violently and has a large difference in magnitude between the two amplitudes (here $1:10^3$). Furthermore, as shown in \figref{fig21}(c) and (d), the high-frequency nature of $A^\varepsilon(x)$ leads to the fact that the solution of the equation and its derivatives have equally strong localized features, which requires that the proposed method needs to be able to capture the localized high-frequency information well.	To capture such oscillations and localized high-frequency information, traditional methods often require very fine meshes, which leads to significant computational cost. The ordinary deep neural network, when solving such problems, captures the oscillations and localized high-frequency information poorly due to the F-Principle, resulting in a relatively poorer computational accuracy. 
	
		\begin{figure}[thbp!]
		\centering
		\includegraphics[width=1\linewidth]{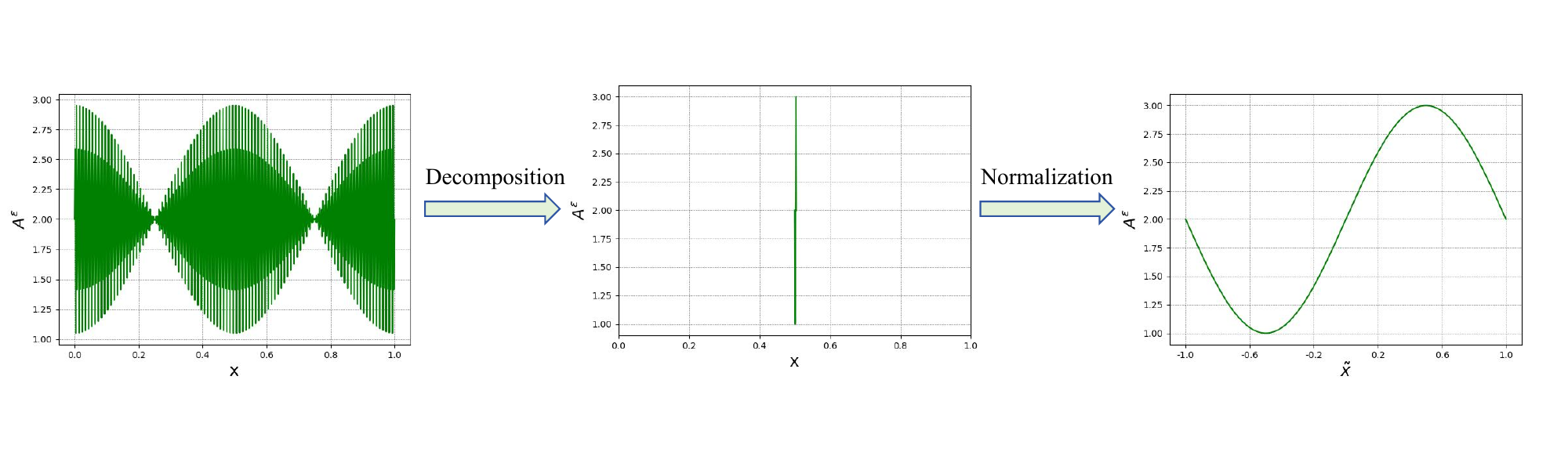} 
		\caption{\footnotesize  Variation of the coefficient $A^\varepsilon(x)$  after domain decomposition and normalization.}
		\label{fig31}
	\end{figure}

	\textbf{Domain decomposition and normalization}. To accommodate small scale ratio, we construct many local solutions by dividing the division of actual computational domain $ \Omega $ into multiple non-overlapping subdomains. By domain decomposition and normalization , the original high-frequency data can be transformed into a low-frequency space. \figref{fig31} demonstrates that the coefficients $A^\varepsilon(x)$ change from high-frequency oscillations to low-frequency in each subdomain after domain decomposition and normalization. Then we approximate the sub-solution by using a local radial basis function neural network, which naturally possesses the local approximation property and makes it possible to better capture the local properties of the solution.
	
	\textbf{Weak formulation and randomized neural network coefficients}.  For the multiscale elliptic equations \eqref{1.1}, we consider  $\boldsymbol{x} \in \Omega \subset \mathbb{R}^n, n=1,2$, and $ f(\boldsymbol{x})$, $g(\boldsymbol{x})$, $u^\varepsilon(\boldsymbol{x}) \in \mathbb{R} $. Since the coefficients $ A^\varepsilon(\boldsymbol{x}) $ and solution $ u^\varepsilon(\boldsymbol{x}) $ are usually functions of $ \boldsymbol{x}/\varepsilon $, when $ \varepsilon $ is small, they have higher frequencies, large order of magnitude differences after derivation, and drastic gradient changes. Solving directly for the strong formulation may pose some difficulties, so we will consider the weak solution and transfer the derivation of $ A^\varepsilon(\boldsymbol{x}) $  and $ \triangledown u^\varepsilon(\boldsymbol{x}) $ to the derivation of the test function $ v(\boldsymbol{x}) $.
	
	Under proper conditions, the weak formulation of the above problems \eqref{1.1} is: Find $ u\in H^1_g(\Omega) $ such that 
	\begin{equation}\label{3.1}
	\int_\Omega (A^\varepsilon\triangledown u) \cdot \triangledown  v dx = \int_\Omega f v dx  \qquad \forall  v \in  H^1_0(\Omega).
	\end{equation} 
	Here, $ H^1_g(\Omega) = \{ u \in H^1(\Omega) : u=g  \text{ on } \partial \Omega   \} $.

	To numerically solve the problem \eqref{3.1}, the trial function space in RRNN method  will be chosen as the randomized radial basis function space $ U_R(\Omega) \subset H^1_g(\Omega) $, while the test function space $ V_h(\Omega) \subset  H^1_0(\Omega)$ is chosen as a finite-dimensional function space, such as finite element spaces.  However, due to the high-frequency characteristics of $ u(\boldsymbol{x}) $, the numerical results of \eqref{3.1} are not satisfactory enough.

	 In order that the original high-frequency data can be transformed into a low-frequency space, we decompose the domain into some subdomains, and approximate the solution on each subdomain by using a local randomized radial basis function neural network. First, we give some notation. Let $ \{\mathcal{T}_h\} $ be a usual partition of $ \Omega $ into $ S $ non-overlapping subdomains, $|\mathcal{T}_h| = S $. Let $ \Gamma $ be the union of the boundaries of all the subdomains $ K\in \{\mathcal{T}_h\} $, $  \Gamma^{in} $ be the set of all interior edges, and $  \Gamma^\partial = \Gamma  \backslash \Gamma^{in} $. Let $ K^{+}_i $ and $ K^{-}_i $ be the two neighboring subdomain sharing a common edge $ e $. 

	The RRNN method is referenced to Petrov-Galerkin finite element method to solve the equations \eqref{1.1}:	In each subdomain $ K $, find $ u_{\rho}\in U_R(K) $ such that 
	\begin{equation}\label{3.2}
		\int_K (A^\varepsilon\triangledown u_{\rho}) \cdot \triangledown  v dx  - \int_{\partial K} (A^\varepsilon\triangledown u_{\rho}) \cdot \boldsymbol{n}_K  v ds =\int_K f v dx \qquad \forall v \in V_h(K),  \quad \forall K \in \mathcal{T}_h. 
	\end{equation}
	where $ \boldsymbol{n}_K $ is the unit outer normal vector on $ \partial K $. However, the solution of Equ.\eqref{3.2} is not equivalent to the multiscale elliptic equations \eqref{1.1} because the local problems lack connections with each other and Dirichlet boundary conditions. To address this deficiency, we randomly choose collocation points on $ \mathcal{T} $ and enforce the solution to satisfy the boundary condition or $ C^0 $ and $ C^1 $ continuity conditions.
	
	We denote the set of collocation points on domain boundary by $ \mathcal{B}=\{  \boldsymbol{x}^b_j \in e: e \in \Gamma^{\partial}, j=1,2,...,N_{b}   \} $ and interior boundary by $ \mathcal{C}=\{  \boldsymbol{x}^c_j \in e: e \in \Gamma^{in}, j=1,2,...,N_{c}   \} $. Let $ K^+ $ and $ K^- $ be two neighboring elements sharing a common edge $ e $, and the unit outward normal vectors on $ \partial K^{\pm} $ be $ \boldsymbol{n}^{\pm} = \boldsymbol{n}|_{\partial K^{\pm}} $.

	Then the RRNN method for solving the problems \eqref{1.1} is: Find $ u_{\rho}\in U=\{u \in L^2(\Omega) : u|_K \in U_R(K), \quad \forall K \in \mathcal{T}_h\} $ such that 
	\begin{align}
		\int_K (A^\varepsilon\triangledown u_{\rho}) \cdot \triangledown  v_k dx  - \int_{\partial K} (A^\varepsilon\triangledown u_{\rho}) \cdot \boldsymbol{n}_K  v_k ds &=\int_K f v_k dx \qquad \forall v_k \in V_h,  \quad \forall K \in \mathcal{T}_h,\label{3.3} \\
		\boldsymbol{n}^{+}u_{\rho}(\boldsymbol{x}^c_j)|_{\partial K^+}  + \boldsymbol{n}^{-}u_{\rho}(\boldsymbol{x}^c_j)|_{\partial K^-} &= 0\quad\qquad\qquad \qquad\qquad\qquad\forall \boldsymbol{x}^c_j \in \mathcal{C}, \label{3.4}\\
		\boldsymbol{n}^{+} \cdot \triangledown u_{\rho}(\boldsymbol{x}^c_j)|_{\partial K^+}  + \boldsymbol{n}^{-}  \cdot \triangledown u_{\rho}(\boldsymbol{x}^c_j)|_{\partial K^-} &= 0 \quad\qquad\qquad\qquad\qquad\qquad\forall \boldsymbol{x}^c_j \in \mathcal{C},\label{3.5} \\
		u_{\rho}(\boldsymbol{x}^b_j) &= g(\boldsymbol{x}^b_j) \qquad\quad \quad \quad\qquad\qquad\forall \boldsymbol{x}^b_j \in \mathcal{B}.\label{3.6}
	\end{align}
	Here $ V_h $ is the test function space composed of the following localized test functions, which defined on a subset $ K\subset \Omega $ reads
	\begin{equation}\label{3.7}
		v_k=\begin{cases}
			\bar{v}\ne 0,\quad \text{over }K,\\
			0, \qquad \quad\text{over } K^c,
		\end{cases}\quad K\cup K^c = \Omega.
	\end{equation}
In one-dimensional, the non-vanishing test function $ \bar{v}(x) $ is chosen to be the high-order polynomials $ P_{k+1}(x)-P_{k-1}(x) $, in which $ P_k(x) $  is Legendre polynomial of order $ k $. In two-dimensional, using the tensor product rule, $ \bar{v}(x,y) $ can be chosen to be $ (P_{k_1+1}(x)-P_{k_1-1}(x))(P_{k_2+1}(y)-P_{k_2-1}(y)) $.

 From Equ.\eqref{3.3}-\eqref{3.6},  we can define the component matrices $ \boldsymbol{A} $ and source term vector $ \boldsymbol{b} $ of the final linear systems to be solved:
\begin{equation}\label{3.8}
	\begin{bmatrix}
			\boldsymbol{A}^{\text{e}} \\
				\boldsymbol{A}^{\text{b}}\\
					\boldsymbol{A}^{\text{c}}
	\end{bmatrix} \boldsymbol{w}
	=\begin{bmatrix}
	\boldsymbol{b}^{\text{e}} \\
	\boldsymbol{b}^{\text{b}}\\
	\boldsymbol{b}^{\text{c}}
\end{bmatrix},
\end{equation}
	where $ \boldsymbol{A}^{\text{e}} $ is a $ SQ\times SJ $ matrix, $ \boldsymbol{b}^{\text{e}}  $ is a $ SQ \times 1 $ vector. Here Q represents the number of test functions $ v_k $ of per subdomain. $ \boldsymbol{A}^{\text{b}} $ is a $ N_{b}\times SJ $ matrix, $ \boldsymbol{b}^{\text{b}}  $ is a $ N_{b} \times 1 $ vector. $ \boldsymbol{A}^{\text{c}} $ is a $ N_{c}\times SJ $ matrix, $ \boldsymbol{b}^{\text{c}}  $ is a $ N_{c} \times 1 $ vector. $ \boldsymbol{w} $ is a $ SJ\times 1 $ vector of unknown weight variables. We abbreviate Equ.\eqref{3.7} as
	\begin{equation}\label{3.9}
		\boldsymbol{A} \boldsymbol{w} = \boldsymbol{b}.
	\end{equation}
	
	To solve the multiscale elliptic equations \eqref{1.1} by the RRNN method, we need to find an optimal parameter $ \boldsymbol{w} $ of the RBFNN to  minimize the error of  Equ.\eqref{3.8}, which is a standard linear least-squares problem, and the  best value of $ \boldsymbol{w}$ is $\hat{\boldsymbol{w}}=\text{argmin}\|\boldsymbol{A} \boldsymbol{w} - \boldsymbol{b}\|. $ 
	
	The algorithm for that the RRNN method solves the multiscale elliptic equations \eqref{1.1} is summarized as follows:
	\begin{algorithm} \label{algorithm1}
		\caption{The RRNN method solves the multiscale elliptic equations \eqref{1.1}}
		\LinesNumbered
		\begin{small}
		\KwIn{ test points $\boldsymbol{x}_{test} $; Parameters  $ N_b $, $ N_c $, $ \beta $, $ S $, $ J $, $Q$.} 
		\KwOut{$ \boldsymbol{w} $; the solution $ u(\boldsymbol{x}_{test}) $.}
		Initialize radial basis function network architecture. \\
		Randomly choose the parameter $ \sigma  $ and $ \boldsymbol{c} $ from a uniform distribution as Equ.\eqref{2.8} and fix them.   \\
		Decompose the domain to $ S $ subdomains and choose $ Q $ test functions $  v_k\in V_h $ as Equ.\eqref{3.7} in each subdomain. \\
		Assemble the linear system $ \boldsymbol{A}^{\text{e}} \boldsymbol{w} =  	\boldsymbol{b}^{\text{e}}  $ according to Equ.\eqref{3.3}. \\
		Randomly take points $ \{\boldsymbol{x}_j^b\}_{j=1}^{N_b} $ and assemble the linear system $ \boldsymbol{A}^{\text{b}} \boldsymbol{w} =  	\boldsymbol{b}^{\text{b}}  $ .\\
		Randomly take points $ \{\boldsymbol{x}_j^c\}_{j=1}^{N_c} $ and assemble the linear system $ \boldsymbol{A}^{\text{c}} \boldsymbol{w} =  	\boldsymbol{b}^{\text{c}}  $ .\\
		Solve linear systems Equ.\eqref{3.8} by linear least-squares method and obtain $ \boldsymbol{w} $.\\
		Compute the solution $ u(\boldsymbol{x}_{test}) $ by Equ.\eqref{2.3}.
			\end{small}
	\end{algorithm}
	
	\section{Numerical examples}	

In this section, we provide a series  of multiscale elliptic numerical examples to test the RRNN method. These examples include one-dimensional and two-dimensional problems, two-scale and three-scale problems, and Poisson-Boltzmann problems.

To evaluate the performance of our proposed method, we have conducted comparative analyses with several prevalent neural network approaches. This includes methods predicated on the gradient descent method and the least squares method. Specifically, we have chosen to compare our method with the Physics-Informed Neural Network (PINN)\cite{pinns} and the Deep Ritz Method (DRM)\cite{drm}, both of which are extensively utilized for the resolution of partial differential equations. Additionally, we have considered the Multiscale Deep Neural Network (MscaleDNN)\cite{MscaleDNN} and the Sparse Radial Basis Function Neural Network (SRBFNN)\cite{srbfnn}, which are tailored for addressing multiscale elliptic problems. Furthermore, the Local Extreme Learning Machine (LocELM) method\cite{locelm}, which operates on principles of domain decomposition and the least squares method, is also included in our comparative analysis. The primary focus of our evaluation is on two critical metrics: computational accuracy and computational efficiency.

We use the following relative maximum error (max error) and relative root mean square error (rms error) to measure computational accuracy:
\begin{equation}\label{4.1}
	\text{max error} = \frac{\|u^N-u\|_{L_{\infty} }}{\| u \|_{L_{\infty} }},\qquad \text{rms error} = \frac{\|u^N-u\|_{L_2}}{\| u \|_{L_2}},
\end{equation}
where $ u^N $ is the approximate solution of the RRNN method and $ u $ is the reference solution. If the numerical example has an exact solution, $ u $ is given directly,  otherwise we use the finite difference method(FDM) to obtain a numerical reference solution $ u $. The grid step size of the FDM is chosen to be 0.0001 for the one-dimensional cases and 0.0005 for the two-dimensional cases. 

Computational cost denotes the aggregate training time for the entire neural network. For the RRNN method and LocELM, this metric encapsulates three distinct phases: pre-processing, optimization, and testing times. The pre-processing phase involves the generation of the coefficient matrix $\boldsymbol{A} $ and the right-hand side source term $ \boldsymbol{b}$ for the least squares problem at hand. This phase encompasses the computational demands associated with deriving the output functions of the network and their respective derivatives, which are efficiently computed leveraging the "\textbf{torch.autograd}" module in Python. The optimization phase pertains to the temporal expenditure on resolving the linear least squares problem, facilitated by the \textbf{scipy.linalg.lstsq} routine. The testing phase,  measures the duration required to ascertain the model's predictive accuracy across a spectrum of test points. For other benchmarked methods, computational cost is defined by the training interval, spanning from the inception to the culmination of the optimization training loop under a prescribed number of epochs and iterations. All timing data	is collected with the "\textbf{time}" module in Python.

For PINN, DRM, MscaleDNN and SRBFNN, they are all trained using Adam\cite{Kingma}. The parameters settings can be found in Appendix A. All the codes for this paper are implemented in Python.

\subsection{One dimensional} 

\noindent \textbf{Example 4.1.} Firstly, we consider the periodic case that $ A^\varepsilon(x) :$ 
\begin{equation}\label{4.2}
	A^\varepsilon(x) = (2 + \cos(\frac{2\pi x}{\varepsilon}))^{-1},
\end{equation}
and the exact solution is given by 
\begin{equation}\label{4.3}
	u(x)=x-x^2+\varepsilon(\frac{1}{4\pi}\sin(\frac{2\pi x}{\varepsilon})  - \frac{1}{2\pi}x\sin(\frac{2\pi x}{\varepsilon})) - \varepsilon^2(\frac{1}{4\pi^2}\cos(\frac{2\pi x}{\varepsilon})  + \frac{1}{4\pi^2} )  
\end{equation}
with boundary condition $  u(0) = u(1) = 0 $ and source term $ f(x)=1 .$ In this study, we aim to evaluate the accuracy and computational efficiency of the RRNN method across various scale ratio $\varepsilon$. 

The primary parameters of the RRNN model under investigation include: the number of subdomains ($ S $), the number of test functions per subdomain ($ Q $), the number of neurons per subdomain ($ J $), and the upper bound of the random shape coefficient ($\beta$). The impact of these different model parameters on the simulation outcomes will be explored, thereby identifying the optimal parameters for subsequent examples. The results can be found in Appendix B. 

For all one-dimensional cases, we will utilize 80 Gauss–Lobatto quadrature points and corresponding weights in performing the integral within each subdomain. \tabref{tab1} presents the experimental findings for RRNN in Example 4.1, showcasing the maximum error, root mean square (rms) error, and training time for varying values of scale ratio $ \varepsilon $.

\begin{center} \begin{minipage}[t]{0.5\linewidth}
		\captionof{table}{\newline\footnotesize   Results of RRNN for Example 4.1 when $ \beta=2, J=100 $ and $ Q=20. $ }	
		\label{tab1}
		\small
		\begin{tabular}{llllc}
			\hline 
			$ \varepsilon $&$ S $ &max error & rms error& \makecell{Training time\\(seconds)}\\ 
			\hline 
			0.5 &5 & 1.48E-13 & 1.25E-13  & 1.74\\
			0.1 &10 & 4.31E-12 &4.17E-12 &2.05\\
			0.05&20 &7.11E-12  &6.49E-12  &2.72 \\
			0.01&50 & 9.87E-11 & 9.03E-11 &4.75\\
			0.005&100  &1.58E-10  &1.57E-10 & 8.67\\
			\hline 	
		\end{tabular}
	\end{minipage}
\end{center}

\begin{figure}[htbp!]
	\centering
	\footnotesize
	\begin{tabular}{@{\extracolsep{\fill}}c@{}c@{\extracolsep{\fill}} @{\extracolsep{\fill}}c }
		\includegraphics[width=0.33\linewidth]{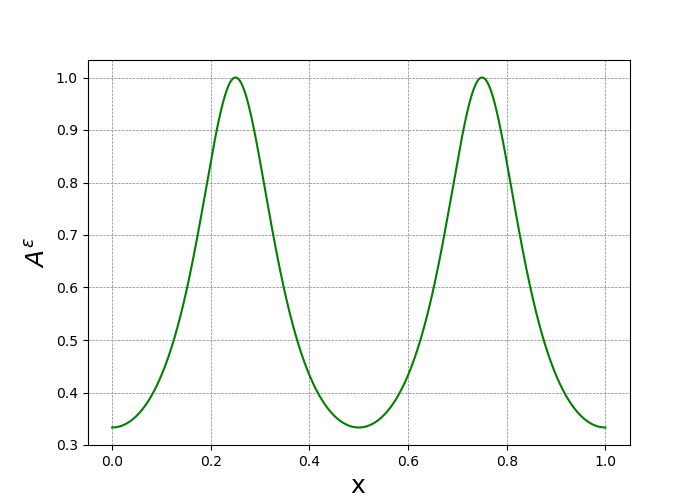} &
		\includegraphics[width=0.33\linewidth]{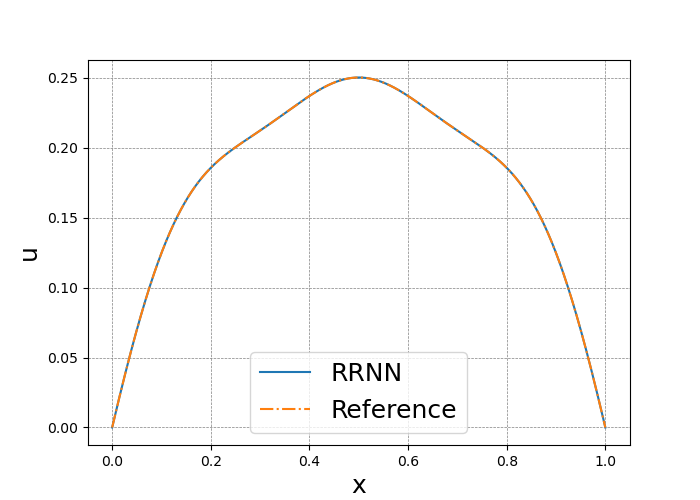}&
		\includegraphics[width=0.33\linewidth]{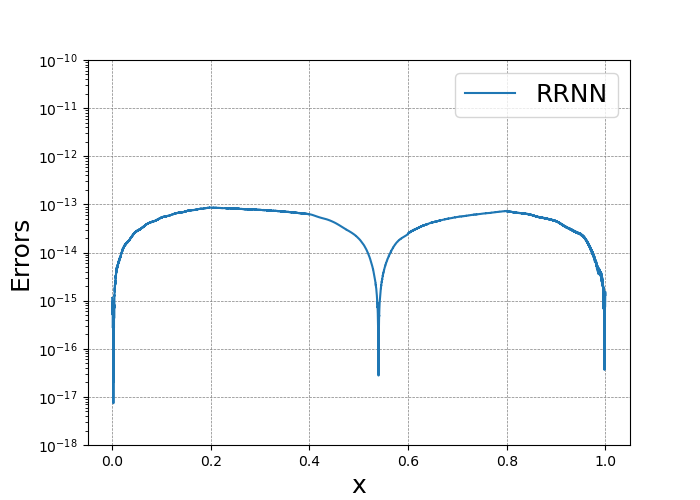} \\
		(a) coefficient $A^\varepsilon(x)$  & (b) solutions&	(c) point-wise error \\
		\includegraphics[width=0.33\linewidth]{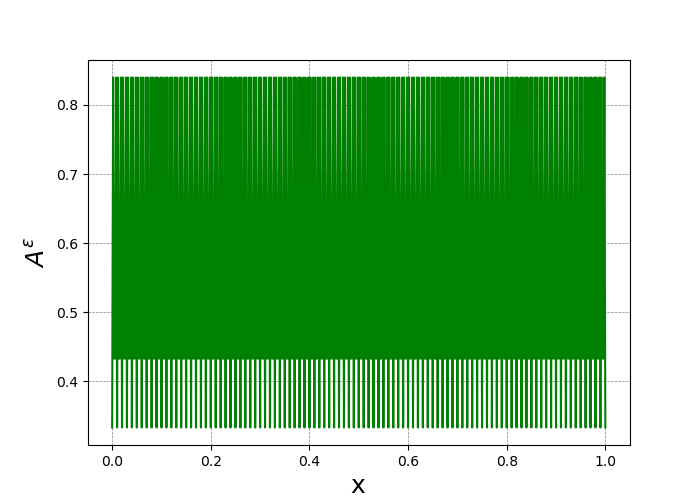} &
		\includegraphics[width=0.33\linewidth]{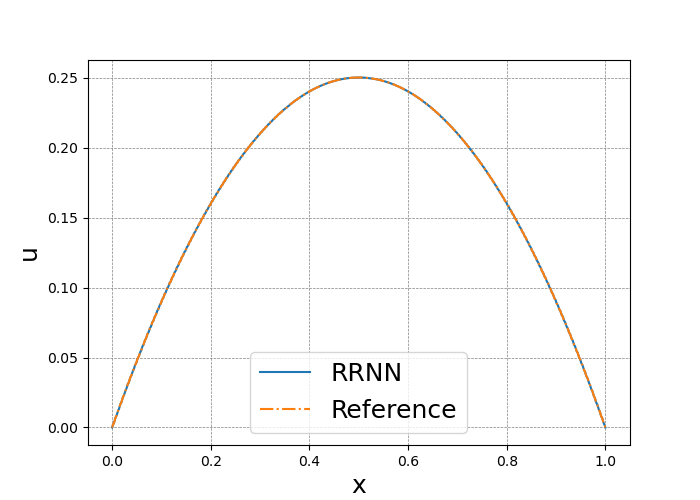}&	
		\includegraphics[width=0.33\linewidth]{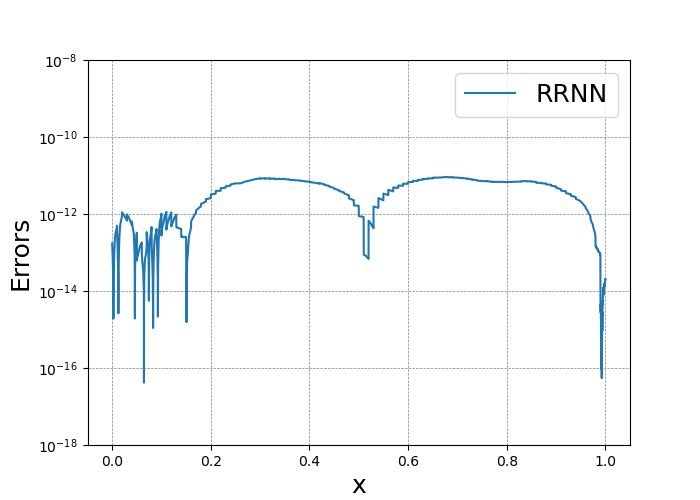}\\	
		(d) coefficient $A^\varepsilon(x)$  & (e) solutions&	(f) point-wise error\\
	\end{tabular}
	\caption{\footnotesize The coefficient $A^\varepsilon(x)$, numerical solution obtained by RRNN and reference solution, and absolute point-wise error for Example 4.1: (a)(b)(c) for $ \varepsilon= 0.5 $, (d)(e)(f) for $ \varepsilon = 0.005$. }
	\label{fig2}
\end{figure}

\figref{fig2} depicts the coefficient $A^\varepsilon(x)$, numerical solutions and corresponding error profiles yielded by the RRNN when $ \varepsilon $ is set to 0.5 and 0.005. The error profiles are defined as the absolute discrepancies between the RRNN solution and the exact solution as given by Equ.\eqref{4.3}. \figref{fig2}(a) and (d) illustrate the values of the coefficient $A^\varepsilon(x)$ on the interval $x\in[0,1]$, which exhibit significant high-frequency characteristics as $\varepsilon$ decreases, which makes solving the equations significantly more difficult. Notably, for this example, the RRNN method is capable of achieving an accuracy on the order of 1E-13, and it maintains an exemplary accuracy level of 1E-10 even when $ \varepsilon $ is as small as 0.005. Furthermore, the method demonstrates a commendable computational efficiency, with the entire training process concluding within a span of 10 seconds.

\noindent \textbf{Example 4.2.} Consider the double-scale case that $ A^\varepsilon(x) :$ 
\begin{equation}\label{4.4}
	A^\varepsilon(x) = 2 + \sin(\frac{2\pi x}{\varepsilon})\cos(2\pi x),
\end{equation}
with $ f(x)=1,g(x)=1. $ Since the exact solution is unknown, we will obtain  a reference solution using the finite difference method(FDM) with a grid step size of 0.0001. 
In this example, RRNN method will be compared with SRBFNN, MscaleDNN, DRM, PINN and LocELM.

For the RRNN method, we employ up to order $ Q=20 $ Legendre polynomials as test functions and $ J=50 $ neurons in each subdomain. The shape coefficient $ \sigma $ is randomly chosen from a uniform distribution of $ \mathbb{U}([0,5]) $. The domain will be divided into $ S= 1/\varepsilon $ subdomain (especially, $ S=5 $ for $ \varepsilon=0.5 $).

\begin{center} \begin{minipage}[t]{0.98\linewidth}
		\captionof{table}{\newline\footnotesize  Computational accuracy comparison between RRNN and deep learning method based on gradient descent in terms of max/rms error for Example 4.2. }	
		\label{tab2}
		\scriptsize
		\centering
		\begin{tabular}{l|lllll|lllll}  
			\hline
			&\multicolumn{5}{c|}{ max error } & \multicolumn{5}{c}{rms error } \\
			\hline 
			$ \varepsilon $& 0.5&0.1& 0.05 & 0.01&0.005 &0.5&0.1& 0.05 & 0.01&0.005 \\
			\hline
			\textbf{RRNN} & 5.82E-8 &3.64E-7&1.22E-6&6.70E-6&1.36E-5 &2.88E-8&1.72E-7&7.06E-7&3.74E-6&7.55E-6 \\	
			
			SRBFNN &1.01E-4 &	5.40E-5&	6.25E-5	&1.84E-4&	2.06E-4  & 2.86E-5&	1.30E-5&	1.85E-5&	5.04E-5&	5.21E-5 \\		
			
			MscaleDNN  &2.44E-3&	2.00E-3&	3.14E-3&	2.87E-3&	2.59E-3&	9.76E-4&	1.16E-3&	2.25E-3&	2.10E-3&	1.77E-3 \\
			
			DRM &	5.97E-3&	6.21E-3&	3.25E-3&	3.45E-3&	3.06E-3&	2.49E-3&	3.43E-3&	2.23E-3&	2.70E-3&	2.22E-3 \\
			
			PINN 	&4.11E-3&	7.77E-3&	7.47E-3&	7.48E-3&	9.31E-3&	1.89E-3&	5.86E-3&	5.72E-3&	5.71E-3&	7.37E-3  \\
			\hline 
		\end{tabular}
	\end{minipage}
\end{center}

\begin{center} \begin{minipage}[t]{0.6\linewidth}
		\captionof{table}{\newline\footnotesize  Computational efficiency comparison between RRNN and deep learning methods based on gradient descent method in terms of training time for Example 4.2. }	
		\label{tab3}	
		\centering	
		\small
		\begin{tabular}{l|lllll}  
			\hline 
			$ \varepsilon $& 0.5&0.1& 0.05 & 0.01&0.005 \\
			\hline
			\textbf{RRNN}&  1.75 & 	1.99 & 	2.72 & 	8.13 & 	17.97 \\	
			
			SRBFNN & 257.77 &	270.14 &	257.81 &	254.15&	261.72  \\		
			
			MscaleDNN  &1310.69 &	1445.23 &	1282.06 &	1255.61 &	1356.51   \\
			
			DRM & 1025.14 &	1026.33 &	1006.53 &	1008.78 &	1060.15   \\
			PINN & 1063.49 & 	1125.37 &	1253.86 &	1177.55 &	1125.20  \\
			\hline 
		\end{tabular}	
	\end{minipage}
\end{center}	

A comparative analysis of the RRNN method against the SRBFNN, MscaleDNN, DRM, and PINN method with respect to maximum/rms error and training time, respectively, is presented in \tabref{tab2} and \tabref{tab3}. In terms of computational accuracy, \tabref{tab2} demonstrates that the RRNN method outperforms the other methods across all scale ratio $\varepsilon$, particularly excelling when the parameter $ \varepsilon $ is large. Notably, the RRNN method achieves a reduction in error by 2-3 orders of magnitude at most scale ratio. While the RRNN method delivers commendable results in accuracy, its most significant advantage lies in computational efficiency. As depicted in \tabref{tab3}, the training time for the RRNN method is on the order of seconds, whereas training times for the MscaleDNN, DRM, and PINN exceed 1000 seconds, and the SRBFNN requires approximately 250 seconds. These findings suggest that the RRNN method is superior to these neural network method based on gradient decent in both approximation accuracy and computational cost.

\begin{center} \begin{minipage}[t]{0.78\linewidth}
		\captionof{table}{\newline\footnotesize Computational accuracy and efficiency comparison between RRNN and neural network method based on least squares method in terms of max/rms error and training time for Example 4.2 when scale ratio $\varepsilon$ is small. }	
		\label{tab4}	
		\centering	
		\small
		\begin{tabular}{lllllllc}  
			\hline 
			$ \varepsilon $ & Method & $ Q $ & $ N $&$ M $  & max error& rms error &\makecell{Training time\\(seconds)} \\
			\hline 
			0.05& \textbf{RRNN}   & 20 & 440 & 1000 & 1.22E-6  & 7.06E-7& 2.72 \\
			& LocELM & 20 & 440 & 1000 & 2.69E-3 & 1.83E-3 & 3.29 \\
			&        & 50 & 1020 & 1000 & 2.57E-6  &  1.28E-6&  3.21\\
			&        & 100 &2040 & 1000 & 2.46E-6 & 6.07E-7& 3.59 \\
			\hline
			0.01& \textbf{RRNN }  & 20 & 2200 & 5000 & 6.70E-6  & 3.74E-6& 1.75 \\
			& LocELM & 20 & 2200 & 5000 &4.62E-1	&3.37E-1	& 5.81 \\
			&        & 50 & 5100 & 5000 &2.78E-3  &  1.99E-3&  10.15\\
			&        & 100 &10200 & 5000 & 1.74E-4 & 1.20E-4  & 10.46 \\
			\hline
			0.005& \textbf{RRNN}   & 20 & 4400 & 10000 & 1.36E-5 &7.55E-6& 17.97   \\
			& LocELM & 20 & 4400 & 10000 & 8.49E-1 & 5.49E-1 & 11.96  \\
			&        & 50 & 10200 & 10000 & 5.05E-3 & 3.43E-3& 52.88\\
			&        & 100 &20400 & 10000 & 6.03E-4 & 3.74E-4  & 54.06 \\
			\hline	
			0.002	& \textbf{RRNN}   & 20 & 11000 & 25000 & 3.45E-5&	2.54E-5&	77.96  \\
			& LocELM & 20 & 11000 & 25000 & 1.05E+0 & 8.59E-1 & 82.27  \\
			&        & 50 & 26000 & 25000 & 1.81E-2   &  1.34E-2&  688.10\\
			&        & 100 &51000 & 25000& 3.75E-2 & 2.75E-2 & 725.12  \\
			\hline
		\end{tabular}	
	\end{minipage}
\end{center}	

The comparative results, as observed in \tabref{tab2} and \tabref{tab3}, indicate that the RRNN method not only has a shorter training time but also offers higher accuracy when compared to these gradient descent-based methods. This advantage is well-understood, as it reflects the general benefits of the least squares method over gradient descent methods. To further elucidate the merits of the RRNN in addressing linear multiscale elliptic problems, we have conducted a comparison with the LocELM method, which also employs the least squares method and domain decomposition. In this comparison, we have maintained parity in the number of subdomains and neurons per subdomain between the LocELM and the RRNN method. For consistency, we denote the number of collocation points used in each subdomain of the LocELM by $Q$. \tabref{tab4} presents the comparative results of RRNN and LocELM in terms of maximum/rms error and training time for Example 4.2, across a range of $\varepsilon$ values from 0.5 to 0.002, which reflects a situation where multiscale effects are relatively strong.  Both methods share the common step of forming a linear system, which is subsequently solved using the least squares approach. In \tabref{tab4}, $N$ and $M$ represent the dimensions of the matrix—the number of rows and columns, respectively—of the linear system to be solved, providing an additional metric that can reflect computational cost.

\begin{figure}[htbp!]
	\centering
	\footnotesize
	\begin{tabular}{@{\extracolsep{\fill}}c@{}c@{\extracolsep{\fill}} @{\extracolsep{\fill}}c }
		\includegraphics[width=0.33\linewidth]{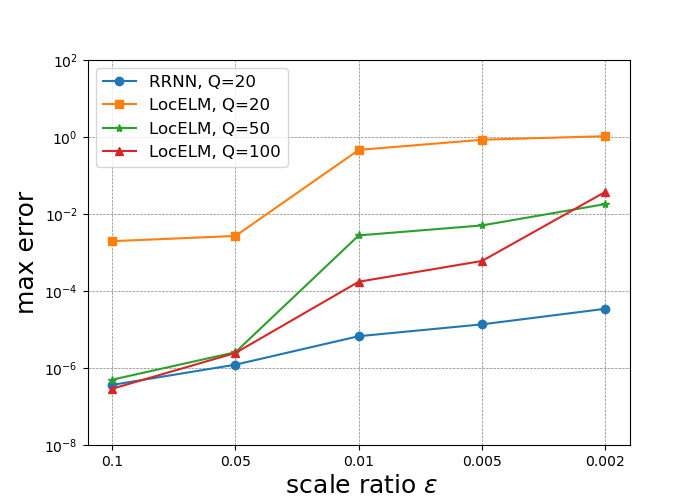} &
		\includegraphics[width=0.33\linewidth]{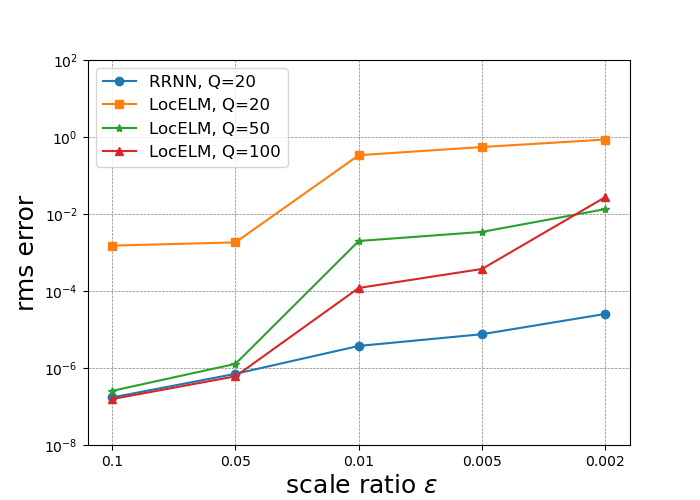}&
		\includegraphics[width=0.33\linewidth]{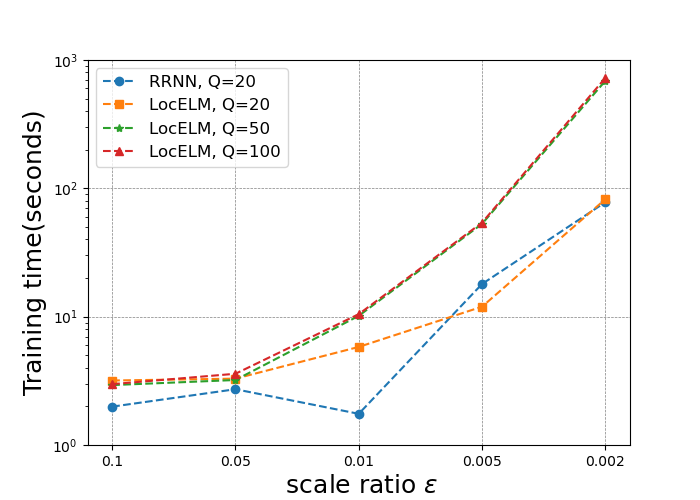}\\	
		(a) max error  & (b) rms error  & (c) training time\\
	\end{tabular}
	\caption{\footnotesize Comparison between RRNN and LocELM with different $Q$ for Example 4.2. }
	\label{fig511}
\end{figure}

The findings presented in \tabref{tab4} and \figref{fig511} reveal that for any given scale ratio $ \varepsilon $, the RRNN method consistently achieves higher accuracy than the LocELM when the dimensions of the matrix to be solved ($N \times M$) are equivalent. When the matrix size is not taken into account but only the precision,  it is observed that for larger values of $\varepsilon$ (e.g., $\varepsilon > 0.01$), both methods exhibit closely matched accuracy and training times, with no significant advantage for the RRNN method. This observation can be attributed to the fact that when $ \varepsilon $ is substantial, the local high-frequency characteristics of the equation are not pronounced, and the benefits of the least squares method become more evident. Conversely, as $ \varepsilon $ diminishes, the LocELM necessitates an increased number of collocation points to attain the desired level of accuracy. This requirement leads to the formulation of a larger-scale matrix, consequently escalating the training time. In stark contrast, the RRNN method maintains a consistent and optimal performance by utilizing a fixed number of collocation points, Q = 20, thereby achieving superior accuracy in a more expedited training time frame. This superiority can be ascribed to the RRNN method's enhanced ability to capture high-frequency features when $ \varepsilon $ is small. The LocELM, which employs the hyperbolic tangent (tanh) activation function and requires the computation of $ A^\varepsilon(\boldsymbol{x}) $ and $ \triangledown u^\varepsilon(\boldsymbol{x}) $, often struggles to encapsulate local high frequency information, thereby limiting its accuracy. In contrast, the RRNN method effectively addresses this limitation.

\noindent \textbf{Example 4.3.} Consider the three-scale problem that $ A^\varepsilon(x) :$ 
\begin{equation}\label{4.5}
	A^\varepsilon(x) = (2 + \cos(\frac{2\pi x}{\varepsilon_1}))(2 + \cos(\frac{2\pi x}{\varepsilon_2})),
\end{equation}
and we take $ f(x)=1 $ and $ g(x)=1. $  Similar to Example 4.2, the reference solution will be generated by the finite difference method(FDM) with a grid step size of 0.0001.

\begin{figure}[htbp!]
	\centering
	\footnotesize
	\begin{tabular}{@{\extracolsep{\fill}}c@{}c@{\extracolsep{\fill}} @{\extracolsep{\fill}}c }
		\includegraphics[width=0.33\linewidth]{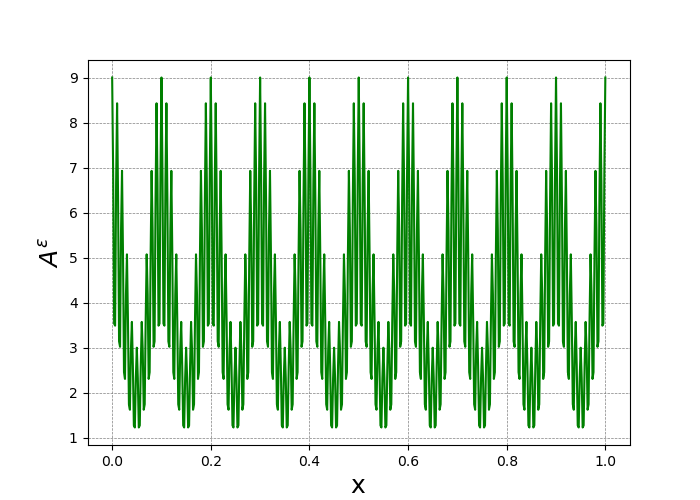} &
		\includegraphics[width=0.33\linewidth]{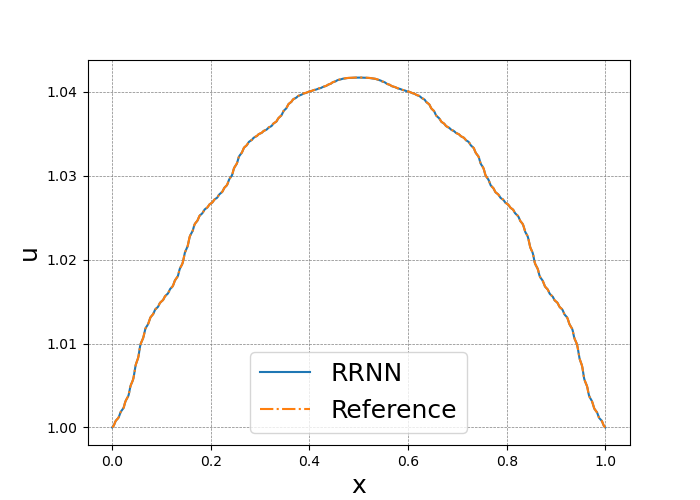}&
		\includegraphics[width=0.33\linewidth]{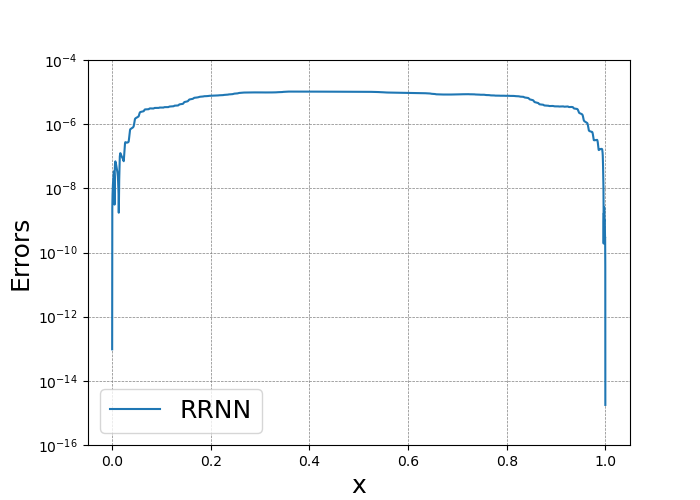}  \\
		(a) coefficient $A^\varepsilon(x)$  & (b) solutions&	(c) point-wise error  \\	
		\includegraphics[width=0.33\linewidth]{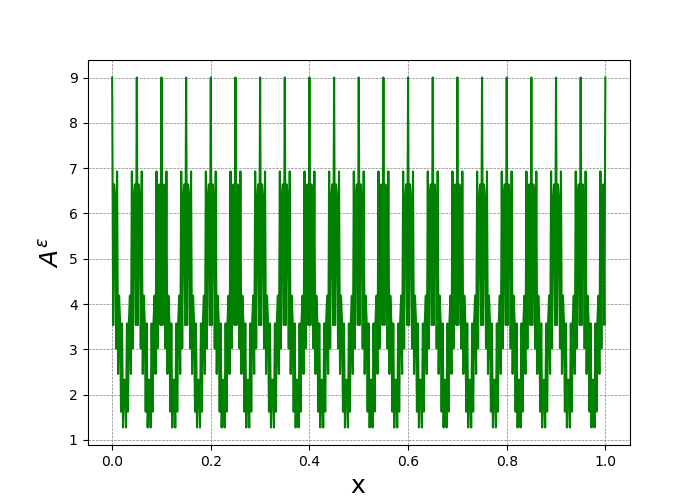} &
		\includegraphics[width=0.33\linewidth]{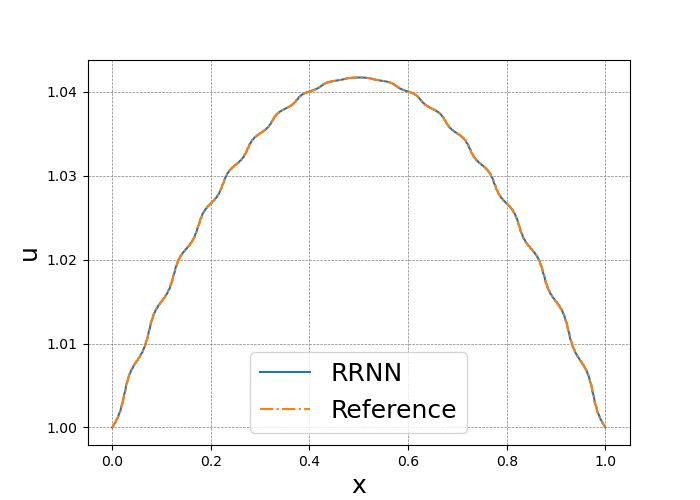}&		
		\includegraphics[width=0.33\linewidth]{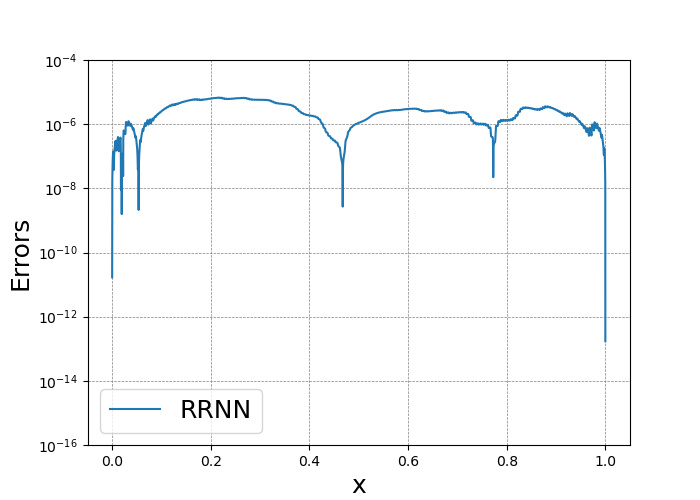}\\	
		(d) coefficient $A^\varepsilon(x)$  & (e) solutions  &	(f) point-wise error \\
	\end{tabular}
	\caption{\footnotesize The coefficient $A^\varepsilon(x)$, numerical solution obtained by RRNN and reference solution, and absolute point-wise error for Example 4.3: (a)(b)(c) for $ \varepsilon_1 = 0.1, \varepsilon_2 = 0.01 $, (d)(e)(f) for $ \varepsilon_1 = 0.05, \varepsilon_2 = 0.005 $. }
	\label{fig51}
\end{figure}

In Example 4.3, we address a example involving a three-scale coefficient as defined by Equ.\eqref{4.5}, with two distinct scale ratio $ \varepsilon_1  $ and $ \varepsilon_2 $. As shown in \figref{fig51}(a) and (d), it oscillates violently and has more complex localized features and variations than the previous two examples , which introduces a heightened challenge in solving the equation. Within the RRNN method simulation framework, we have utilized 100 and 200 uniform subdomains $ (S = 100, 200) $ for scenarios where $ \varepsilon_1, \varepsilon_2 $ are set to $ 0.1, 0.01 $ and $ 0.05, 0.005 $, respectively. The experiment employs $ Q = 20 $ Legendre polynomials as test functions and $ J = 50 $ neurons per subdomain. Additionally, we have set $ \beta = 5  $ to generate the random shape coefficient $ \sigma $. The outcomes of these two distinct experimental setups are detailed in \figref{fig51}, \tabref{tab5} and \tabref{tab71}.

\figref{fig51} presents the  coefficient $A^\varepsilon(x)$ and the numerical solutions procured by the RRNN alongside the reference solutions for Example 4.3. Specifically, \figref{fig51}(a) and (d) illustrates the coefficients of violent oscillations in two cases. \figref{fig51}(b) and (e) depict the RRNN solution profiles for the cases where $ \varepsilon_1 = 0.1, \varepsilon_2 =0.01 $ and $ \varepsilon_1 = 0.05, \varepsilon_2 = 0.005 $, respectively. \figref{fig51}(c) and (f) illustrate the corresponding absolute error profiles. The results indicate that the errors associated with the RRNN method in solving this three-scale example can be within the order of 1E-5 to 1E-6.

\begin{center} \begin{minipage}[t]{0.85\linewidth}
		\captionof{table}{\newline\footnotesize Computational accuracy and efficiency comparison between RRNN and deep learning methods based on gradient descent method  for Example 4.3.}	
		\label{tab5}	
		\centering	
		\small
		\begin{tabular}{l|l|lllll}  
			\hline 
			$ (\varepsilon_1,\varepsilon_2) $&&\textbf{ RRNN}  &SRBFNN  & MscaleDNN &DRM  &PINN \\
			\hline
			(0.1, 0.01)	&max error  &4.79E-6		&	5.80E-5	&	4.36E-3	&	8.69E-3		&	1.57E-2	 \\		
			&rms error &	2.55E-6		&	1.96E-5	&	3.06E-3	&	6.31E-3		&	1.15E-2	\\

			&Training time &	6.52 		&	272.84 	&	1330.88 	&	1011.05 	 &		1150.36	\\
			
			\hline 
			
			(0.05, 0.005)	&max error  &2.73E-5		&	8.02E-5	&	4.87E-3	&	9.07E-3		&	1.63E-2		 \\		
			&rms error &	2.04E-5		&	2.72E-5	&	3.39E-3	&6.44E-3		&	1.21E-2		\\

			&Training time &	13.99 		&	283.69	&	1263.87	&	1009.55		&	1173.45	 	\\
			
			\hline 
		\end{tabular}
	\end{minipage}
\end{center}

\tabref{tab5} offers a computational accuracy and efficiency comparison between RRNN and deep learning methods based on gradient descent method  for Example 4.3, under the conditions that $ \varepsilon_1, \varepsilon_2=0.1,0.01 $ and $0.05,0.005  $. It is observed that the RRNN frequently outperforms MscaleDNN, DRM, and PINN in terms of accuracy, and is comparable to SRBFNN. Furthermore, our method demonstrates a marked improvement in computational efficiency and reduces the training time by almost two orders of magnitude.

\begin{center} \begin{minipage}[t]{0.84\linewidth}
		\captionof{table}{\newline\footnotesize Computational accuracy and efficiency comparison between RRNN and neural network method based on least squares method for Example 4.3.}		
		\centering		
		\label{tab71}	
		\small
		\begin{tabular}{l|l|lllll}  
			\hline 
			&	& \textbf{RRNN} &\multicolumn{4}{c}{LocELM} \\
			$ (\varepsilon_1,\varepsilon_2) $&	$ (Q , J) $& $ (50, 20) $ &$(50, 20) $  & $ (50, 100)  $ & $(100, 20)  $ &$ (100, 100)  $ \\
			\hline
			
			(0.1, 0.01)	&max error  &4.79E-6	&	7.49E-2	&	1.31E-2	&1.10E+0	&	5.72E-6		 \\		
			&	rms error &	2.55E-6	&	5.12E-2 	&	9.34E-3	&9.68E-1	&	4.23E-6	\\

			&	Training time &	6.52 	&	4.47  	&	9.24	&	6.94 	&	43.64 	\\
			
			\hline 
			
			(0.05, 0.005)	&	max error  &2.73E-5	&1.00E+0	&	3.47E-2	&	1.10E+0	&	2.24E-5			 \\		
			&	rms error &	2.04E-5	&	9.78E-1	&	2.52E-2	&	9.86E-1	&	1.66E-5			\\

			&	Training time &	13.99 	&	9.60 	&	44.77 	&	15.04 	&	273.20 		\\
			
			\hline 
		\end{tabular}
	\end{minipage}
\end{center}

The above results demonstrate the advantages of the  RRNN  method over the neural network method based gradient descent. For a fairer comparison, we will compare it with the least squares-based method. \tabref{tab71} presents a comparative evaluation of the LocELM with varying values of $ Q $ and $ J $ against the RRNN method for Example 4.3. The number of subdomains is consistent with that used in the RRNN method. The results underscore that the RRNN method achieves significantly better computational accuracy for a comparable computational cost (i.e., same matrix size). While augmenting the values of $ Q $ and $ J $ can enhance the accuracy of LocELM, it also results in a pronounced escalation in training time. Conversely, the RRNN method successfully addresses the problem within a single-digit time, given the same level of desired accuracy.

\subsection{Two dimensional} 
\noindent \textbf{Example 4.4.} Consider the case:
\begin{equation}\label{4.6}
	A^\varepsilon(x,y) =\frac{1}{4+\cos(\frac{2\pi(x^2+y^2)}{\varepsilon})},
\end{equation}	
with smooth source function $ f(x,y)=-(x^2+y^2) $. The exact solution is 
\begin{equation}\label{4.7}
	u(x,y) =\frac{1}{4}(x^2+y^2)^2+\frac{\varepsilon}{16\pi}(x^2+y^2)\sin(\frac{2\pi(x^2+y^2)}{\varepsilon}) + \frac{\varepsilon^2}{32\pi^2}\cos(\frac{2\pi(x^2+y^2)}{\varepsilon})
\end{equation}
and boundary conditions $ g(x,y) $ can be defined from the exact solution.

In all two-dimensional examples, We employ Legendre polynomials in each
direction $ x $ and $ y $, and perform the integral in each element by employing $ 10\times10 $ Gauss-Lobatto quadrature
points using tensor product rule. In each subdomain, test function $ {v}_k(x,y) $ will be chosen to be ${v}_k(x,y)= (P_{k_1+1}(x)-P_{k_1-1}(x))(P_{k_2+1}(y)-P_{k_2-1}(y))(k_1=1,2,...,Q;k_2=1,2,...,Q) $, where $ P_k(x)  $ is  Legendre polynomials of order $ k $, where uses $ Q $ test functions in each direction $ x $ and $ y $. The random collocation points on domain boundary $ N_b  $ and on interior boundary $ N_c $ will be fixed as $ N_b = 4SN_{bper} $ and $ 2S(S-1)N_{cper} $, where $ N_{bper} $ and $ N_{cper} $ represent the number of random collocation points in each edge of per subdomain and will take values of $  N_{bper} = N_{cper} = 10  $. 
\begin{figure}[htbp!]
	\centering
	\footnotesize
	\begin{tabular}{@{\extracolsep{\fill}}c@{}c@{\extracolsep{\fill}} @{\extracolsep{\fill}}c }
		\includegraphics[width=0.33\linewidth]{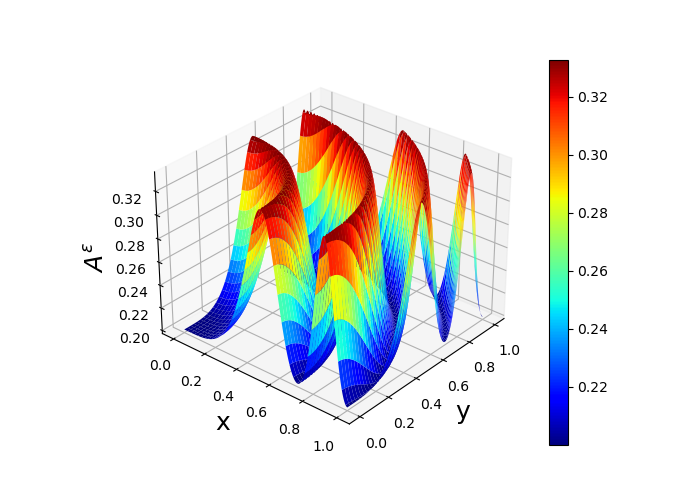} &
		\includegraphics[width=0.33\linewidth]{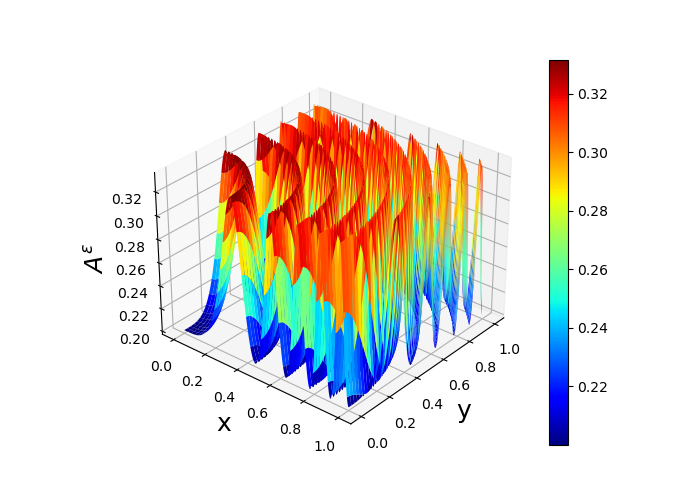}&
		\includegraphics[width=0.33\linewidth]{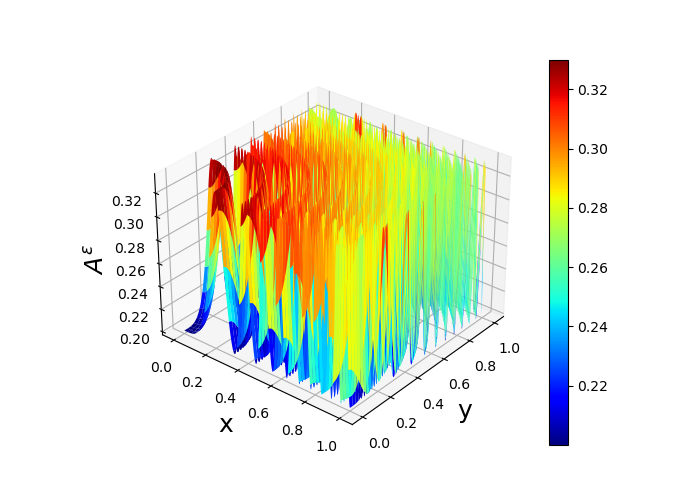}\\	
		(a) coefficient $A^\varepsilon(x,y)(\varepsilon=0.5)$ & (b) coefficient $A^\varepsilon(x,y)(\varepsilon=0.2)$  & (c) coefficient $A^\varepsilon(x,y)(\varepsilon=0.1)$\\
		\includegraphics[width=0.33\linewidth]{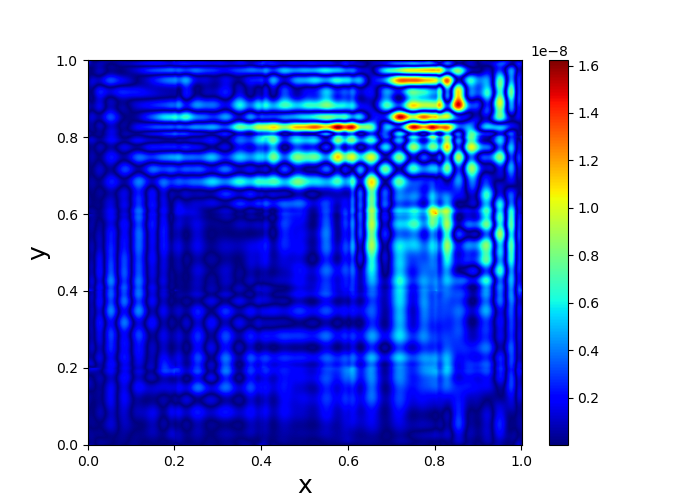} &
		\includegraphics[width=0.33\linewidth]{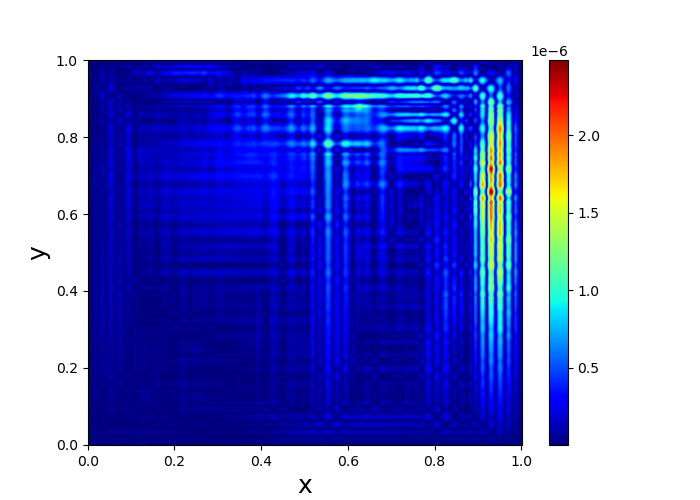}&
		\includegraphics[width=0.33\linewidth]{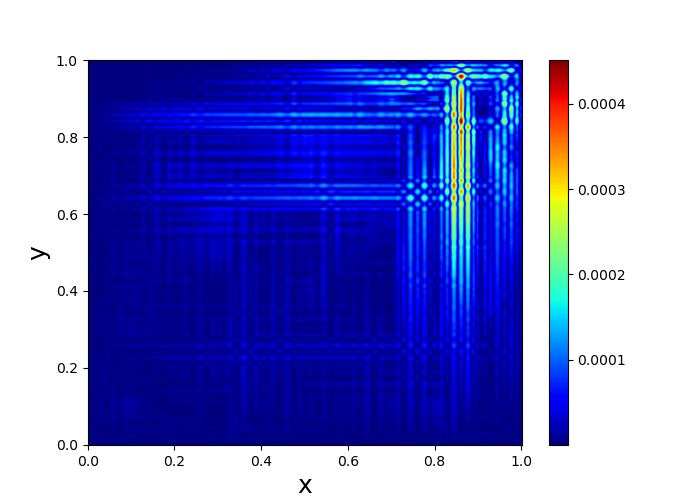}\\	
		(d) point-wise error of  $\varepsilon=0.5$   & (e) point-wise error of  $\varepsilon=0.2$  & (f) point-wise error of $\varepsilon=0.1$\\
	\end{tabular}
	\caption{\footnotesize The coefficient $A^\varepsilon(x,y)$ and absolute point-wise errors of the RRNN method for Example 4.4. }
	\label{fig6}
\end{figure}

In this example, the RRNN method will be compared with LocELM and the results are shown in \tabref{tab9}. In RRNN method, the upper bound of random shape coefficient $ \beta = 1$. In both methods, the number of neurons will be fixed on $ J=200 $, and the domain is composited into $ 5\times5, 8\times8,10\times10$ for $ \varepsilon=0.5,0.2,0.1 $, respectively. 

 \figref{fig6} presents the  coefficient $A^\varepsilon(x,y)$ and the absolute point-wise errors for varying scale ratio within the context of Example 4.4. It is observable that the absolute point-wise errors are predominantly within the orders of 1E-8, 1E-6, and 1E-4 for $ \varepsilon = 0.5, 0.2, 0.1 $, respectively. These observations substantiate the capability of our proposed method to effectively approximate solutions for two-dimensional multiscale elliptic equations with a coefficient of violent oscillation by a high precision and efficiency.

\tabref{tab9} delineates a comparative analysis between the RRNN method and  the LocELM when applied to Example 4.4, with both methods employing $ Q = 10  $ test functions or collocation points in each direction. The number of subdomains and the number of neurons per subdomain for the LocELM are aligned with those utilized in the RRNN method. When the matrices of the linear systems solved by both methods are of equivalent size — implying that the computational cost is close — the RRNN method consistently demonstrates superior accuracy across all scales. These findings suggest that the RRNN method possesses a more favorable potential for addressing linear multiscale elliptic problems.

\begin{center} \begin{minipage}[t]{0.55\linewidth}
		\captionof{table}{\newline\footnotesize Computational accuracy and efficiency comparison of RRNN and LocELM in terms of max/rms error and training time for Example 4.4  when both are $ Q=10 $. }	
		\label{tab9}	
		\centering	
		\small
		\begin{tabular}{llllc}  
			\hline 	
			$ \varepsilon $ & Method  & max error& rms error &\makecell{Training time\\(seconds)} \\
			\hline 
			0.5	& \textbf{RRNN}    &1.62E-8&	1.26E-8&	16.04 \\
			& LocELM & 1.52E-5	&2.38E-5&	14.39  \\
			\hline
			0.2	& \textbf{RRNN}   &  2.48E-6&	1.34E-6&	70.17  \\
			& LocELM &  1.18E-4&	1.42E-4&	70.32  \\
			\hline
			0.1& \textbf{RRNN}   &  4.51E-4&	2.19E-4&	203.10 \\
			& LocELM  &3.90E-3&	5.15E-3&	216.60 \\
			\hline
		\end{tabular}	
	\end{minipage}
\end{center}

\noindent \textbf{Example 4.5.} Consider a coefficient with double-scale
\begin{equation}\label{4.8}
	A^\varepsilon(x,y) = \frac{1.5+\sin(\frac{2\pi x}{\varepsilon})}{1.5+\sin(\frac{2\pi y}{\varepsilon})} + \frac{1.5+\sin(\frac{2\pi y}{\varepsilon})}{1.5+\cos(\frac{2\pi x}{\varepsilon})} + \sin(4x^2y^2) +1
\end{equation}
with $ f(x)=-10,g(x)=0. $ The exact solution cannot be stated and a reference solution will be gotten by FDM with a step size 0.0005. In this case, we will compare our method with PINN, DRM, MscaleDNN, SRBFNN and LocELM for values of $ \varepsilon $ equal to 0.5, 0.2 and 0.1.  The parameters are consistent with Example 4.4 except for $ \beta = 3 $ and the test functions in each direction $ Q = 9 $.

The RRNN method addresses the equations by resolving a linear system through the least squares method. The accuracy and training time of this method are contingent upon the characteristics and scale of the matrix involved. Specifically, the scale of the matrix being solved are $ 3425 \times 5000, 8864 \times 12800, $ and $ 13900 \times 20000 $ for $ \varepsilon $ values of 0.5, 0.2, and 0.1, respectively. As $ \varepsilon $ diminishes, the high-frequency components of the solution become increasingly pronounced, which exacerbates the complexity of the equation. Consequently, an incremented number of subdomains are required, leading to larger matrices to be solved. This results in a corresponding increase in training time and a subsequent decrease in accuracy.

\begin{figure}[htbp!]
	\centering
	\footnotesize
	\begin{tabular}{@{\extracolsep{\fill}}c@{}c@{\extracolsep{\fill}} @{\extracolsep{\fill}}c }
		\includegraphics[width=0.33\linewidth]{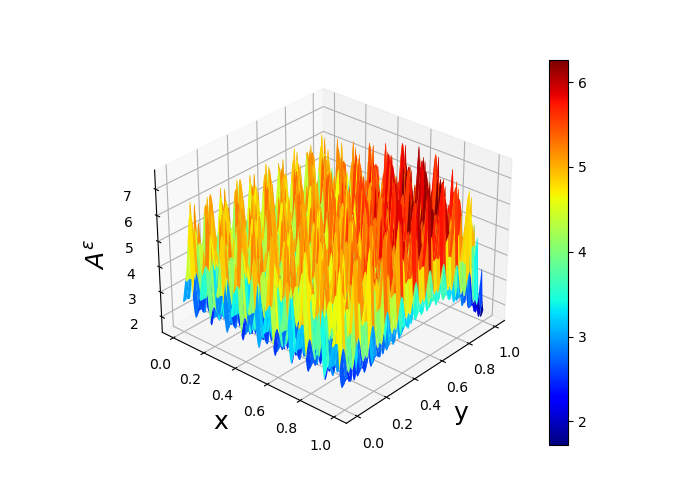} &
		\includegraphics[width=0.33\linewidth]{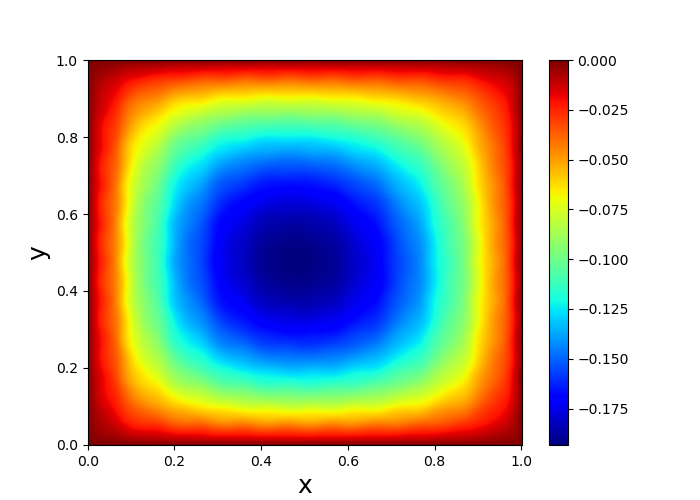}&
		\includegraphics[width=0.33\linewidth]{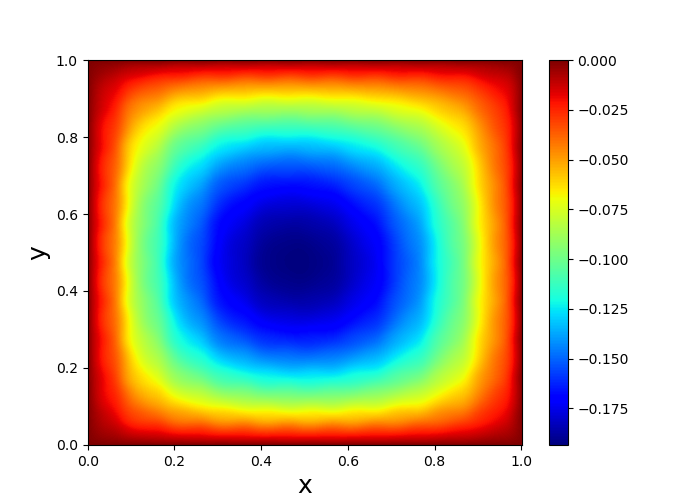}\\	
		(a) coefficient  $A^\varepsilon(x,y)$  & (b) reference solution  & (c) numerical solution of RRNN\\
		\includegraphics[width=0.33\linewidth]{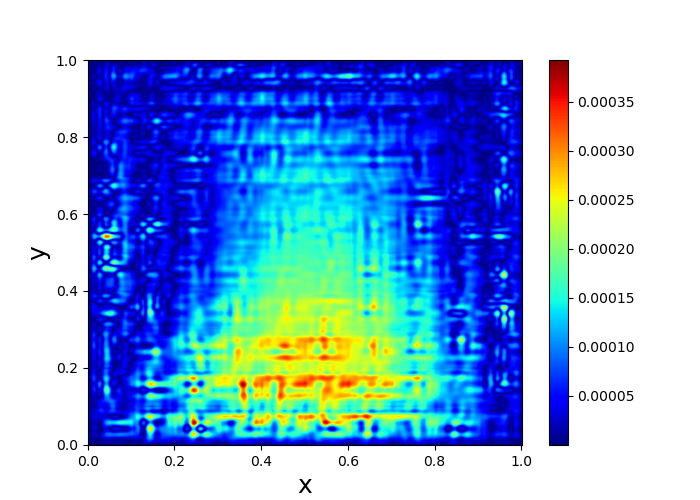} &
		\includegraphics[width=0.33\linewidth]{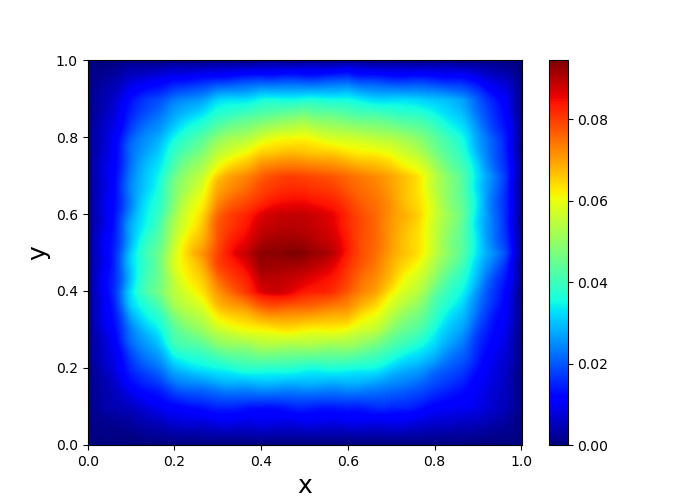}&
		\includegraphics[width=0.33\linewidth]{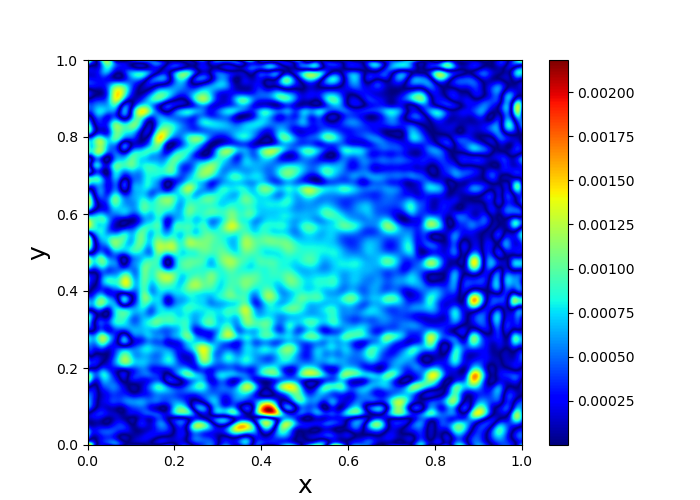}\\	
		(d) point-wise error of RRNN  & (e) point-wise error of LocELM & (f) point-wise error of SRBFNN\\
		\includegraphics[width=0.33\linewidth]{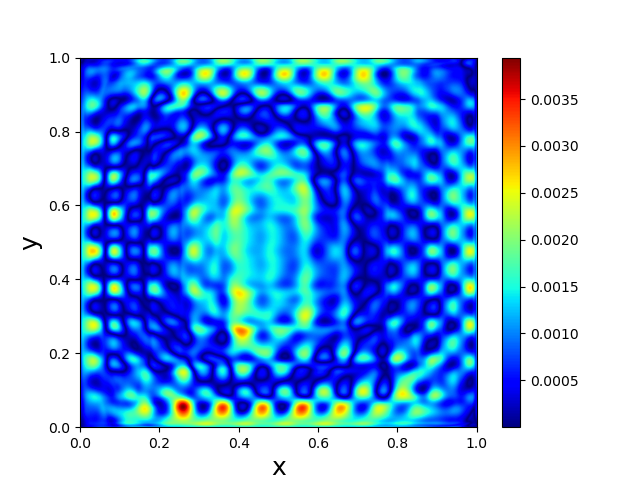} &
		\includegraphics[width=0.33\linewidth]{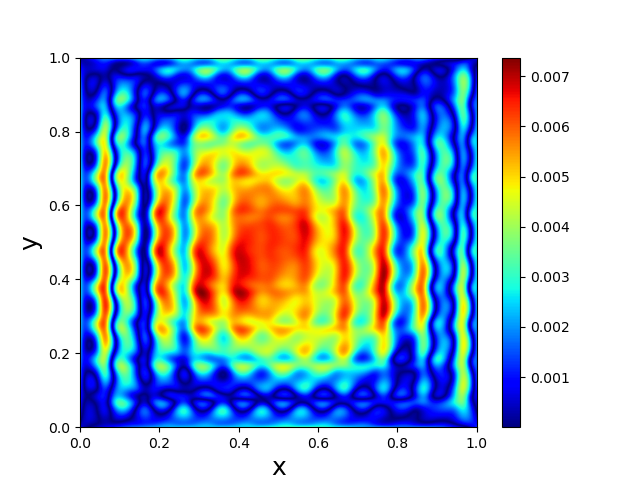}&
		\includegraphics[width=0.33\linewidth]{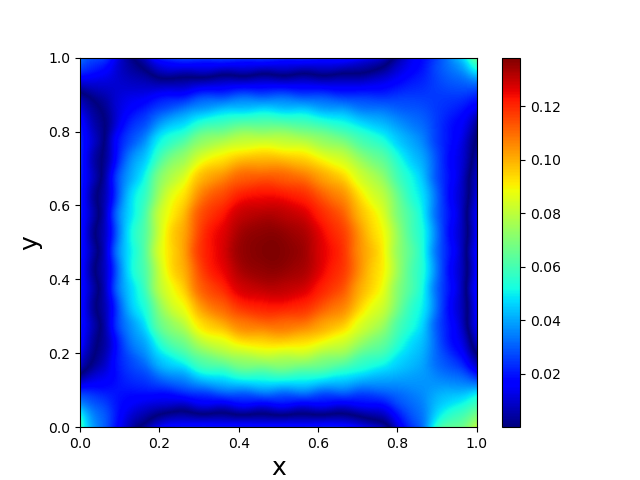}\\	
		(g) point-wise error of MscaleDNN   & (h) point-wise error of DRM   & (i) point-wise error of PINN \\
	\end{tabular}
	\caption{\footnotesize The coefficient $A^\varepsilon(x,y)$, reference solution, numerical solution and  point-wise errors results of Example 4.5 when $\varepsilon=0.1$. }
	\label{fig71}
\end{figure}

In this example, the coefficient $A^\varepsilon(x,y)$ have two different frequency components and is violently oscillating as shown in \figref{fig71}(a) when $\varepsilon = 0.1.$ As shown in \figref{fig71}(e) and (i), the ordinary DNN-based or ELM-based method tend to fail to solve such multiscale problems. Due to the high-frequency nature of the violently oscillating coefficients, a better solution accuracy can only be obtained if targeted improvements are made to the network, as shown in the results presented in \figref{fig71}(d), (f), (g) and (h). Among the solution results of all methods shown in \figref{fig71}, it can be seen that the absolute point-wise error in  \figref{fig71}(d) is the lowest, which is the result of the RRNN method.

\tabref{tab10} presents a more detailed comparative analysis of the performance of the RRNN, LocELM, SRBFNN, MscaleDNN, DRM, and PINN in terms of maximum/rms error and training time for Example 4.5. In this instance, the foundational parameter settings for the LocELM are maintained in alignment with those of the RRNN to ensure an equitable comparison with respect to computational cost. The comparative results demonstrate that the RRNN method achieves superior accuracy, with errors approximately 1-2 orders of magnitude lower than those of the other methods. Furthermore, the RRNN method exhibits significantly greater computational efficiency than the gradient descent-based methods, with training times approximately 1-2 orders of magnitude faster.

\begin{center} \begin{minipage}[t]{0.95\textwidth}
		\captionof{table}{\newline\footnotesize  Computational accuracy and efficiency comparison between  RRNN and some baseline neural network methods  in terms of max/rms error and training time for Example 4.5. }	
		\label{tab10}
		\footnotesize
		\centering
		\begin{tabular}{l|lll|lll|lll}  
			\hline
			&\multicolumn{3}{c|}{ max error } & \multicolumn{3}{c|}{rms error } & \multicolumn{3}{c}{Training time(seconds)}\\
			\hline 
			$ \varepsilon $& 0.5&0.2& 0.1 & 0.5&0.2& 0.1 & 0.5&0.2& 0.1  \\
			\hline
			\textbf{RRNN} &1.42E-3&	6.70E-4&	2.06E-3 &3.30E-4&	3.16E-4&	1.08E-3&	11.76 &	59.51  &	209.00 \\	
			LocELM 	&1.85E-2&	3.57E-2	&3.86E-1&	7.66E-2&	2.59E-2&	3.58E-1&	12.35 &	88.01 &	236.97 
			\\
			SRBFNN &1.86E-2&	1.45E-2&	1.13E-2&	2.05E-3&	3.77E-3&	5.37E-3&	1354.19 &	1227.23 &	1382.60 
			\\		
			
			MscaleDNN  &1.77E-2 &	1.65E-2 &	2.04E-2 &	7.40E-3 &	9.60E-3 &	1.01E-2 &	4614.25  &	4501.84 & 	4656.73 
			\\
			
			DRM &	4.44E-2&	4.90E-2&	3.81E-2&	2.39E-2&	2.76E-2&	3.13E-2&	2829.71 &	2847.55 &	2851.14 
			\\
			
			PINN 	&8.66E-2&	6.92E-1	&7.15E-1&	1.53E-2&	3.13E-1&	6.31E-1&	4762.21 &	4876.50 &	4965.05 
			\\
			\hline 
		\end{tabular}
	\end{minipage}
\end{center}

\noindent \textbf{Example 4.6.} We finally consider the following Poisson-Boltzmann equation with Dirichlet boundary condition,
\begin{subequations}  \label{4.9}
	\begin{numcases}
		\text{-}\text{div}(A^\varepsilon(x,y)\triangledown u^\varepsilon(x,y)) + \kappa(x,y)u^\varepsilon(x,y) = f(x,y), \qquad (x,y) \in \Omega,\label{4.9a} \\
		\qquad \qquad\qquad\qquad\qquad\quad\quad\text{ } \qquad u^\varepsilon(x,y) = g(x,y), \qquad (x,y) \in \partial\Omega,\label{4.9b}
	\end{numcases} 
\end{subequations}
where the inverse Debye-Huckel length of an ionic solvent $ \kappa(x,y)=\pi^2 $ and the dielectric constant $ A^\varepsilon(x,y) =1+0.5\cos(10\pi x)\cos(20\pi y) $. We impose an exact solution $ u(x,y)=\sin(\pi x)\sin(\pi y)+0.05\sin(10\pi x)\sin(20\pi y) $, which determines the corresponding $ g(x,y) $ and $ f(x,y) .$

The Poisson-Boltzmann (PB) equation\cite{caiComputational} plays an important role in many applications including the study of disease and drug design. The PB equation \eqref{4.9} adds an “$ \kappa(x,y)u^\varepsilon(x,y) $” to the multiscale elliptic equations \eqref{1.1}. In order to solve this equation by RRNN, the weak formulation in each subdomain similar to equation \eqref{3.3} will be built first. We multiply the test function $ v_k $ by each term in equation \eqref{4.9a} and integrate in subdomain $ K $ and use a partitioned integral to obtain:
\begin{align}
	\int_K (A^\varepsilon\triangledown u_{\rho}) \cdot \triangledown  v_k dx  - \int_{\partial K} (A^\varepsilon\triangledown u_{\rho}) \cdot \boldsymbol{n}_K  v_k ds   +  \int_{K} \kappa u_{\rho}v_kdx =\int_K f v_k dx \quad \forall v_k \in V_h,  \quad \forall K \in \mathcal{T}_h,\label{4.10}
\end{align}
Combining with equations \ref{3.4}-\eqref{3.6}, the PB equation \eqref{4.9} can be solved by Algorithm 1. To solve this equation, the RRNN method sets the upper bound of random shape coefficient $ \beta = 2$. The number of neurons and test functions in each direction of per subdomain and subdomains will be fixed on $ J=200 $, $ Q=9 $ and $ S=10\times 10 $. 

\begin{figure}[htbp!]
	\centering
	\footnotesize
	\begin{tabular}{@{\extracolsep{\fill}}c@{}c@{\extracolsep{\fill}} }
		\includegraphics[width=0.33\linewidth]{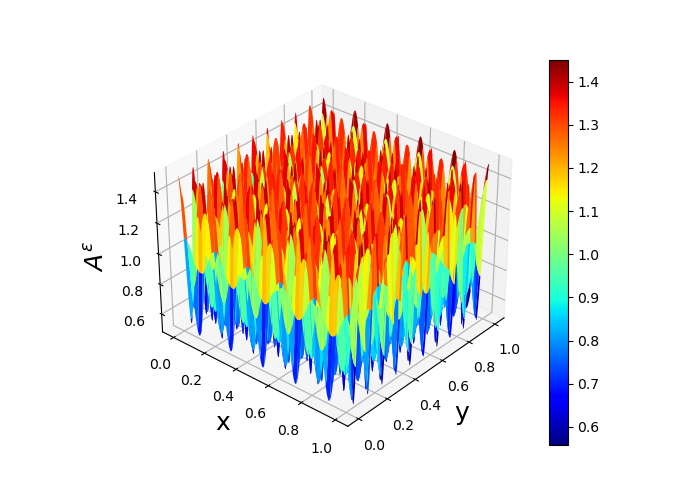} &
		\includegraphics[width=0.33\linewidth]{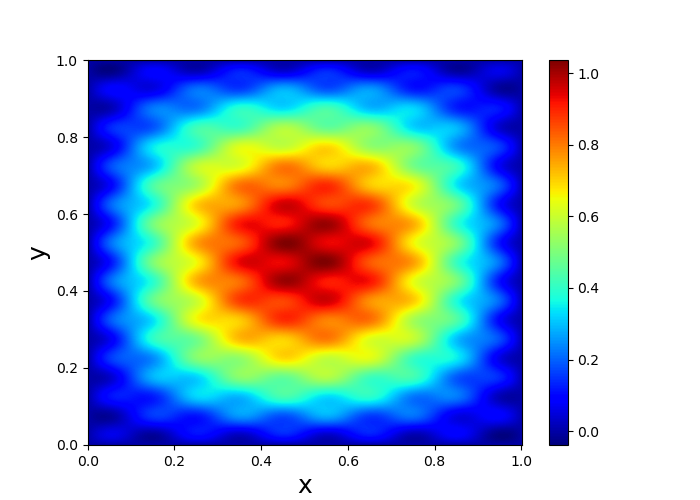} \\
		(a) coefficient $A^\varepsilon(\boldsymbol{x})$ &(b) exact solution of PB equation\\
		\includegraphics[width=0.33\linewidth]{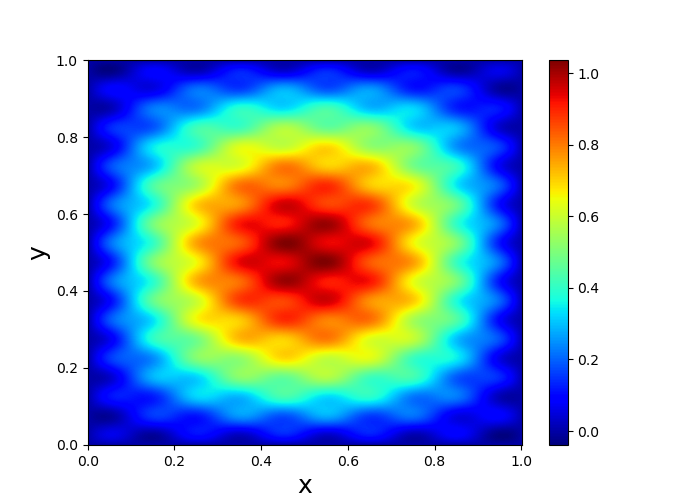} &
		\includegraphics[width=0.33\linewidth]{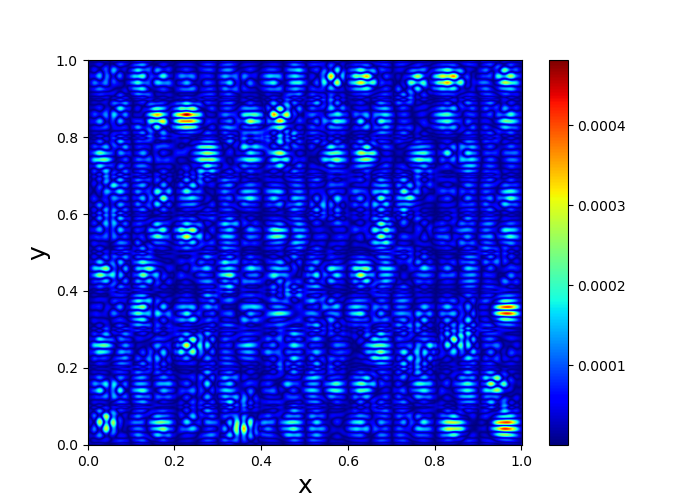}\\
		(c) numerical solution of RRNN&(d) point-wise error of RRNN \\
	\end{tabular}
	\caption{\footnotesize  The coefficient $A^\varepsilon(\boldsymbol{x})$, exact solution, numerical solution and the absolute point-wise errors obtained by RRNN for PB equation in Example 4.6. }
	\label{fig7}
\end{figure}

\figref{fig7}(a) gives the values of the coefficients $A^\varepsilon(x,y)$ of this problem over the solution domain, which exhibits strong oscillations that will increase the difficulty for neural network-based method to address multiscale PDEs\eqref{1.1}. The exact solution  and numerical solution for the Poisson-Boltzmann (PB) equation, as defined by Equ.\eqref{4.9}, are illustrated in Figure \ref{fig7}(b) and \ref{fig7}(c). \figref{fig7}(d) presents the absolute point-wise errors derived from the RRNN method, which shows that the absolute point-wise errors are predominantly on the order of 1E-4 across most regions. This affirms that the RRNN method is capable of effectively solving the PB equation with small scale ratio with high accuracy.

\begin{center} \begin{minipage}[t]{0.93\linewidth}
		\captionof{table}{\newline\footnotesize Computation accuracy and efficiency comparison between  RRNN and some baseline neural network methods for PB equation.}		
		\centering		
		\label{tab11}	
		\begin{tabular}{l|llllll}  
			\hline 
			& \textbf{RRNN}& LocELM&SRBFNN  & MscaleDNN &DRM  &PINN \\
			\hline
			max error & 4.41E-4	& 1.56E-2 	&4.83E-2	&	5.27E-2	&	5.79E-2	&	1.60E+0		 \\		
			rms error &1.43E-4	&1.26E-2		&2.49E-2	 &	3.55E-2&	5.03E-2&	1.09E+0	
			\\	
			training time &192.54&213.12	&1293.48 &	4715.71	&	3116.24 &		4712.07 	
			\\
			
			\hline 
		\end{tabular}
	\end{minipage}
\end{center}

To further substantiate the superiority of the RRNN method, a comparative analysis of the maximum/rms error and training time for several neural network methods is detailed in \tabref{tab11}. The data elucidates that the RRNN method is significantly more precise, offering accuracy four orders of magnitude higher than the PINN methods, and two orders of magnitude higher than that of other methods,. In addition to its superior accuracy, the RRNN method is also markedly more efficient in terms of network training time, which is nearly 5-20 times shorter when compared to the other methods except for LocELM under consideration.

\section{Conclusion} 

Multiscale problems are of paramount importance in the realms of scientific and engineering research, where machine learning algorithms have exhibited promising potential in tackling such complexities. However, challenges pertaining to approximate accuracy and computational cost have often constrained their broader application in multiscale problem-solving. In this paper, we introduce an efficient methodology named RRNN for addressing multiscale elliptic problems. The method adeptly captures the high-frequency information and local characteristics of the solution to multi-scale equations through a combination of regional decomposition, normalization, and the application of a radial basis neural network. Furthermore, it incorporates the concept of extreme learning machines by introducing random coefficients, which significantly curtail the number of training parameters and expedite the training process.

The RRNN method composites the domain and employs a variational formulation in each subdomain. It utilizes the random radial basis function  neural network space approximation as the trail function space, and can flexibly choose the appropriate test function space. The subdomains are interconnected with both $ C^0 $ and $ C^1 $ continuity conditions, combining with  boundary conditions, ultimately culminating in the formulation of a unified linear system concerning the weights of the output layer. This system can be efficiently resolved using the linear least squares method. The extensive examples presented throughout this paper consistently demonstrate the RRNN method's superiority in terms of accuracy and computational efficiency when compared to gradient descent-based methods. Furthermore, in comparison to extreme learning method-based approaches, the RRNN method achieves higher accuracy with equivalent computational cost, with the advantage being particularly pronounced at smaller scale ratios.

We are confident in the potential of this novel method for solving linear multiscale elliptic equations, yet there are several aspects that necessitate further exploration. The first question pertains to the efficient resolution of cases with sufficiently  smaller scale ratios. As the parameter $ \varepsilon $ decreases, our method entails the division of an increasing number of subdomains, resulting in a larger matrix to be solved and consequently, an extended training time.  It is crucial to investigate whether a more optimal approach can be devised to circumvent this challenge. Additionally, it would be beneficial to establish a theoretical framework analogous to finite element analysis to analyze the error associated with the RRNN method.

\subsection*{Appendix A. Model and training setup} 
 In all numerical examples, we consider the computational domain $ \Omega=[0,1]^n,(n=1,2) $. The training details of all models in the numerical experiments are elaborated in the following:
\begin{itemize}
	\item \textbf{PINN:} For one-dimensional examples, it has five linear layers with the activation function tanh. The number of neurons in the hidden layer is 16, 16, 32, 32, 64  when $ \varepsilon $ =0.5, 0.1, 0.05, 0.01, 0.005. We sample 10000 training points and 10000 test points in the domain with a uniform	distribution. The penalty coefficient of the boundary loss function is 20.0, batch-size is 2048, Max-Niter is 3000, and the learning rate is 0.0001. For two-dimensional examples, it has nine linear layers with the activation function tanh. The number of neurons in the hidden layer is 64, 64, 128 when $ \varepsilon $ =0.5, 0.2, 0.1. We sample 500$ \times $500 training points and 1000$ \times $1000 test points in the domain with a uniform	distribution. The penalty coefficient of the boundary loss function is 10.0, batch-size is 2048, Max-Niter is 1000, and the learning rate is 0.001. 
	\item \textbf{DRM:} In one-dimensional and two-dimensional examples, it has six and ten linear layers, respectively, and two residual blocks with the activation function tanh. The penalty coefficients of the boundary loss function are 200.0 and 500.0, respectively. Other parameters are consistent with PINN. 
	\item \textbf{MscaleDNN: } It has eight linear layers with the activation function s2ReLU,
	$$
	\text{s2ReLU}( x ) = \text{sin}  (2\pi x)  * \text{ReLU}( x )   *  \text{ReLU}( 1 - x ). $$	 
	The scale vector is $ [1,2,\cdots,\frac{1}{\varepsilon}]$ and the number of neurons in the hidden layer is (1000, 200, 150, 150, 100, 50, 50) and (1000, 400, 300, 300, 200, 100, 100) for all $ \varepsilon $, respectively, in one-dimensional and two-dimensional examples. In Example 4.5, the learning rate is 0.001, and the penalty coefficient of the boundary loss function is 10000,1000,1000 when $ \varepsilon = 0.5, 0.1, 0.2 $. In Example 4.6, the learning rate is 0.01, and the penalty coefficient of the boundary loss function is 1000. The remaining parameters are consistent with DRM.
	\item \textbf{SRBFNN: } In all one-dimensional examples, the initial learning rate is set to be 0.1 and reduced
	by 1/10 every 300 iterations if it is more than $ 10^{-5} $. We sample 10000 training points and 10000 testing points in the domain with a uniform distribution. tol$_1 $ is  $ 10^{-3} $, tol$_2 $ is  $ 10^{-5} $. The batchsize is 2048, MaxNiter is 3000, SparseNiter is 2000, Checkiter is 100. The hyperparameters $\lambda_1,\lambda_2 $ and $ \lambda_3 $ are initialized to 1.0, 100.0 and 0.001, respectively. The initial number of RBFs is set to 100, 200, 300, 500, 1000, 1500 for six different scales, respectively.  In Example 4.5 and Example 4.6, the initial learning rate is set to be 0.1 and 0.001, respectively, and it is reduced by 1/10 every 30 iterations if it is more than $ 10^{-5} $. We sample 500$ \times $500 training points and 1000$ \times $1000 test points in the domain with a uniform distribution. The batchsize is 2048 in the domain and 512 × 4 on the boundary, MaxNiter is 300, SparseNiter is 250, Checkiter is 10. The hyperparameters tol$_1 $, tol$_2 $, $\lambda_1,\lambda_2 $  are initialized to 0.1, $ 10^{-5} $, 0.01, 20.0, respectively.  $ \lambda_3 $  is initialized to 0.001 and 0.0001 for Example 4.5 and Example 4.6, respectively. The initial number of RBFs is set to 1000, 1000 and 2000 when $ \varepsilon =0.5$, $  0.2 $ and $ 0.1 $, respectively.
	\item \textbf{LocELM: } The network has a single hidden layer with activation function tanh. Domain decomposition and number of neurons will be consistent with the RRNN method. The max magnitude of random coefficients $ R_m $ of LocELM has an influence on the accuracy as discussed in \cite{locelm}, and we have tested and selected a suitable $ R_m $. In the order of LocELM's experiments, $ R_m $ takes the values 3, 5, 1, 3, 3 and 1 in order. 
\end{itemize}  

\subsection*{Appendix B. The influence of diverse parameters}

To ascertain the influence of diverse parameters on the experimental outcomes and to facilitate the selection of optimal parameters for experiments, an investigation into the impact of varying parameters on accuracy and training duration was conducted for Example 4.1.

\begin{figure}[htbp!]
	\centering
	\footnotesize
	\begin{tabular}{@{\extracolsep{\fill}}c@{}c@{\extracolsep{\fill}} }
		\includegraphics[width=0.33\linewidth]{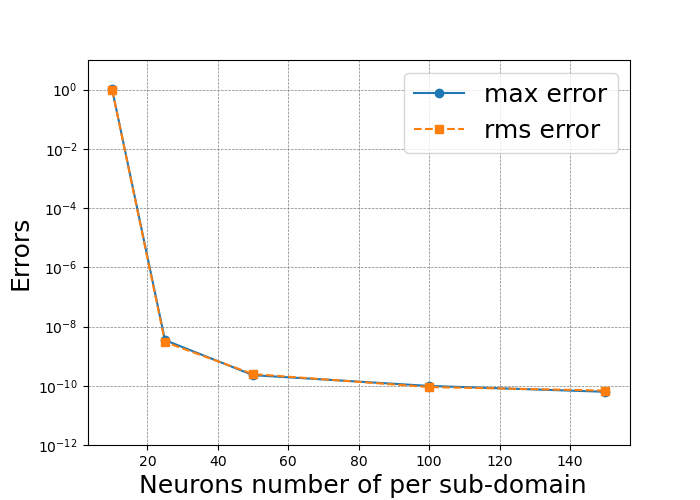} &
		\includegraphics[width=0.33\linewidth]{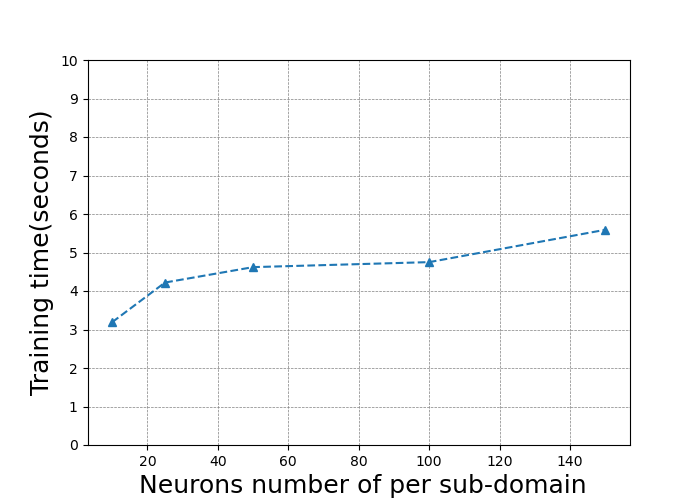}\\
		(a) max error and rms error  & (b) training time \\
		\includegraphics[width=0.33\linewidth]{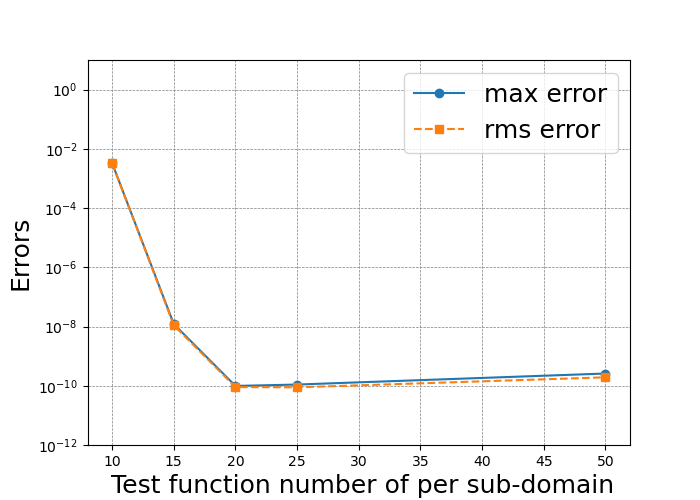} & 
		\includegraphics[width=0.33\linewidth]{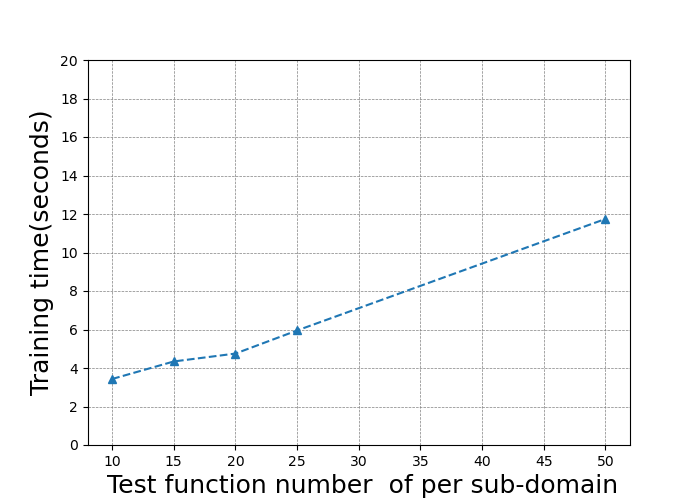}\\	
		(c) max error and rms error& (d) training time\\
	\end{tabular}
	\caption{\footnotesize Effect of the neurons number $ (J) $ and test functions number $ (Q) $ of per subdomain for Example 4.1 :  (a) The max error and rms error in the domain as a function of the neuron number of per subdomain, and (b) the training time; (c) The max error and rms error in the domain as a function of  the test functions number of per subdomain, and (d) the training time.}
	\label{fig3}
\end{figure}

\figref{fig3}(a) and (b) elucidate the influence of the neuron number per subdomain for Example 4.1.  These testes were undertaken with parameters set to $ \varepsilon = 0.01, \beta = 2, Q = 20, S = 50 $. An initial increase in the neuron number per subdomain yields an exponential reduction in numerical error. However, once the number surpasses a threshold (approximately $ J \sim 50 $ for this instance), the rate of error diminishment decelerates. \figref{fig3}(b) suggests that training time ascends marginally in approximate linearity with the neuron count per subdomain, extending from approximately 3.19 seconds for 10 neurons to about 5.59 seconds for 150 neurons.

\figref{fig3}(c) and (d) depict the effects of varying the number of test functions per subdomain. With parameters $ \varepsilon = 0.01, \beta = 2, S = 50, J = 100  $ held constant, they showcase the relative errors and the training time in the domain as a function of the number of test functions per subdomain. \figref{fig3}(c) indicates that augmenting the number of test functions per subdomain initially results in an exponential decrease in numerical error, with the rate of reduction plateauing and potentially reversing slightly once a certain count is exceeded (around $ Q \sim 20 $ in this scenario). \figref{fig3}(d) demonstrates that the training time increases linearly, albeit at a minor rate, with the number of test functions per subdomain.

\begin{figure}[htbp!]
	\centering
	\footnotesize
	\begin{tabular}{@{\extracolsep{\fill}}c@{}c@{\extracolsep{\fill}}  }
		\includegraphics[width=0.33\linewidth]{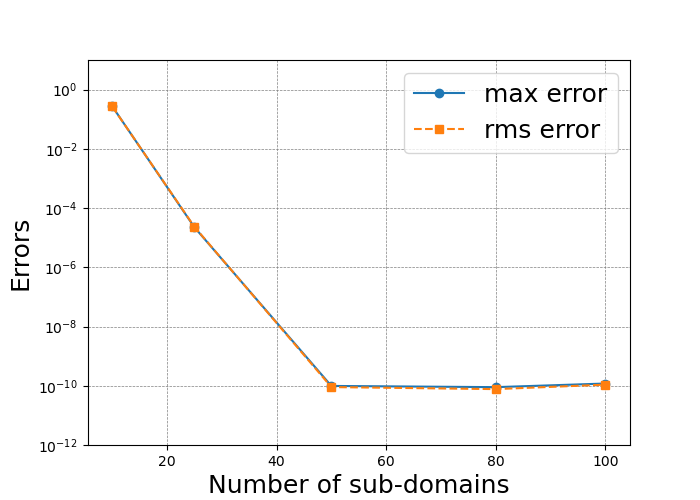} &
		\includegraphics[width=0.33\linewidth]{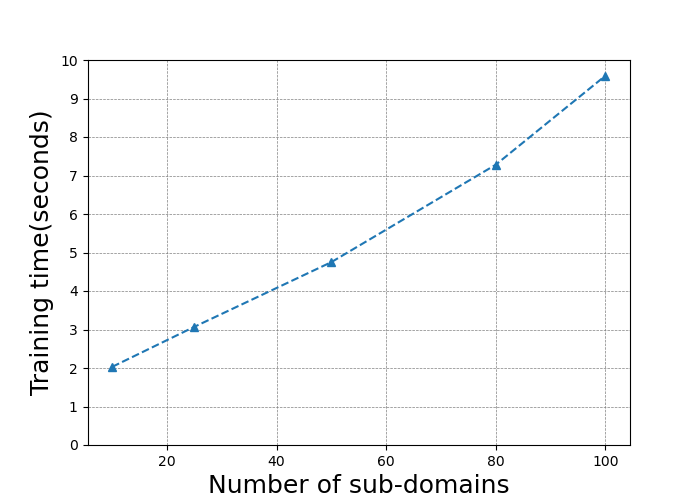}\\
		(a) max error and rms error & (b) training time \\
		\includegraphics[width=0.33\linewidth]{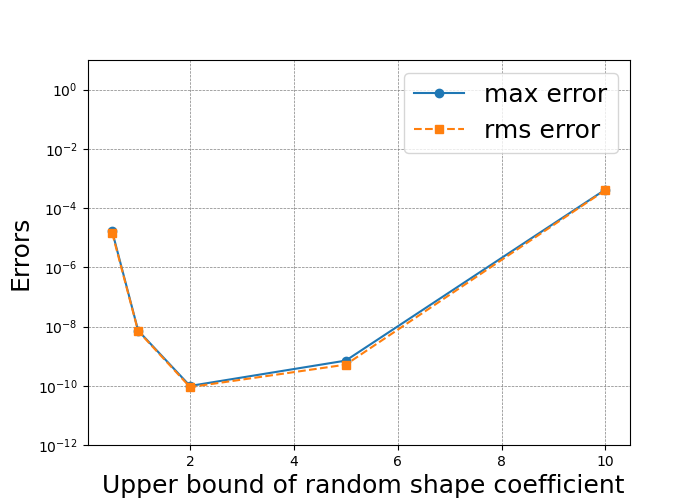} & 
		\includegraphics[width=0.33\linewidth]{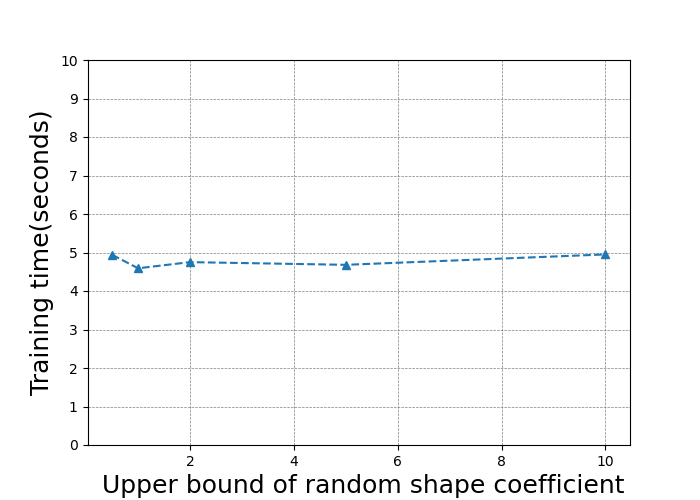}\\	
		(c) max error and rms error & (d) training time\\
	\end{tabular}
	\caption{\footnotesize Effect of the number of subdomains $ (S) $ and the upper bound of random shape coefficient $\sigma (\beta) $ for Example 4.1 : (a) The max error and rms error in the domain as a function of the number of subdomains, and (b) the training time; (c) The max error and rms error	in the domain as a function of upper bound of random shape coefficient $ \sigma $, and (d) the training time.}
	\label{fig4}
\end{figure}

\figref{fig4}(a) and (b) examine the effects of the number of subdomains on the RRNN method simulation, with parameters $ \varepsilon = 0.01, \beta = 2, J = 50, Q = 20 $ maintained. \figref{fig4}(a) portrays the relative maximum error and the relative root mean square error of the RRNN solution across the entire domain. \figref{fig4}(b) delineates the training time for the entire neural network. It is observed that the accuracy of the numerical simulation outcomes procured by the RRNN method generally escalates incrementally with an increasing number of subdomains. Initially, the accuracy improves exponentially with the addition of subdomains; however, after a certain value is attained, further enhancements in accuracy become challenging. Empirical results suggest that optimal accuracy is typically achieved when the number of subdomains equals $ 1/\varepsilon $. Conversely, the total training time progresses approximately linearly with the increment of subdomains, with duration ranging from approximately 2.03 seconds for 10 subdomains to about 9.59 seconds for 100 subdomains.

Ultimately, the investigation reveals that the value of the random shape coefficient $ \sigma $ of the Gaussian Radial Basis Function (GRBF) significantly influences the accuracy of the RRNN method simulation outcomes. As previously discussed in Section 2.1, the shape coefficient of the GRBF within the radial basis function neural network is pre-set to a uniformly random value within the interval $ [0, \beta] $ and remains static throughout the computation. It has been observed that the upper bound $ \beta $ of the random shape coefficient  markedly affects simulation precision. \figref{fig4}(c) and (d) highlight this effect through systematic testing, with other parameters fixed at $ \varepsilon = 0.01, J = 100, Q = 20, S = 50 $, and $ \beta $ varied within the set $ [0.5, 1, 2, 5, 10] $. These figures exhibit the errors and the training time as a function of $ \beta $. The findings indicate that the error tends to deteriorate when $ \beta $ is either substantially large or small, while the training time remains largely unaffected.

\end{document}